\documentclass[10pt]{amsart}
\usepackage{amssymb}
\usepackage{amsbsy}
\usepackage{amscd}
\usepackage{amsmath}
\usepackage{amsthm}

\theoremstyle{plain}
 \newtheorem{thm}{Theorem}[section]
 \newtheorem{lem}[thm]{Lemma}
 \newtheorem{prop}[thm]{Proposition}
 \newtheorem{cor}[thm]{Corollary}
 
\theoremstyle{definition}
 \newtheorem{defn}{Definition}[section]
\theoremstyle{remark}
 \newtheorem{rem}{Remark}[section]
 
 \newtheorem{claim}{Claim}[section]
\def\Bbb{\mathbb}
\def\frak{\mathfrak}
\def\cal{\mathcal}

\newcommand{\Ext}{\operatorname{Ext}}
\newcommand{\Hom}{\operatorname{Hom}}

\newcommand{\im}{\operatorname{im}}

\newcommand{\rk}{\operatorname{rk}}

\newcommand{\NS}{\operatorname{NS}}
\newcommand{\coker}{\operatorname{coker}}
\newcommand{\Pic}{\operatorname{Pic}}
\newcommand{\ch}{\operatorname{ch}}
\newcommand{\td}{\operatorname{td}}

\newcommand{\Hilb}{\operatorname{Hilb}}

\newcommand{\IT}{\operatorname{IT}}
\newcommand{\WIT}{\operatorname{WIT}}

\font\b=cmr10 scaled \magstep5
\def\bigzerou{\smash{\lower1.7ex\hbox{\b 0}}}

\numberwithin{equation}{section}
\begin{document}
\title{
A note on Fourier-Mukai transform
}
\author{K\={o}ta Yoshioka}
 
\address{
Department of mathematics, Faculty of Science, Kobe University,
Kobe, 657, Japan}

\email{yoshioka@math.kobe-u.ac.jp}
 \maketitle

\section{Introduction}\label{sect:intro}

Let $X$ be an abelian or a K3 surface defined over ${\Bbb C}$.
For a smooth projective variety $Z$,
${\bf D}(Z)$ denotes the bounded derived category 
of coherent sheaves on $Z$.
For an abelian surface or a K3 surface $Y$ and an object 
${\cal E} \in {\bf D}(X \times Y)$, an integral functer
\begin{equation}
\begin{matrix}
{\cal F}_{\cal E}: {\bf D}(X)& \to & {\bf D}(Y)\\
\quad x & \mapsto & {\bf R}p_{Y*}(p_X^*(x) \otimes {\cal E})
\end{matrix}
\end{equation}
is called the Fourier-Mukai transform, if 
${\cal F}_{\cal E}$ is an equivalence of categories,
where 
$p_X$ and $p_Y$ are projections from $X \times Y$ to
$X$ and $Y$ respectively.
The Fourier-Mukai transform is a very useful tool for analysing
the moduli spaces of sheaves on $X$.
In order to apply the Fourier-Mukai transform
to an actual problem, it is important to study the problem on
the preservation of stability under the Fourier-Mukai transform. 
We assume that $Y$ is a fine moduli space of sheaves on $X$ and 
${\cal E}$ is the universal family.
In \cite{Y:7}, \cite{Y:9}, we discussed this problem
and showed that the stability is preserved provided
a suitable twisted degree is 0 or 1.
In this note, we show that the Fourier-Mukai transform does not always 
preserve the stability, even for
a $\mu$-stable vector bundle.
We construct two examples (see sect. \ref{sect:counter}):
Assume that $X$ is an abelian surface and $\widehat{X}$
the dual of $X$.
Let ${\cal P}$ be the Poincar\'{e} line bundle on
$X \times \widehat{X}$.
Our first example is constructed for the Fourier-Mukai functor
${\cal F}_{\cal P}$ originally considered by Mukai \cite{Mu:2}. 
We next construct an example for the Fourier-Mukai functor
on a K3 surface.

In section \ref{sect:asymptotic}, we shall provide positive results
on this problem.  
Let $H$ be an ample divisor on $X$.
For a coherent sheaf $E$ on $X$,
$R^i p_{Y*}(p_X^*(E(mH)) \otimes {\cal E})=0$, $i>0$
for $m \gg 0$.
Hence the Fourier-Mukai transform of $E(mH)$, $m \gg 0$ is a sheaf. 
Under some assumptions 
we shall show that
the Fourier-Mukai transform preserves the stability
(cf. Theorem \ref{thm:asymptotic}, Theorem \ref{thm:asymptotic2}).
In \ref{subsect:basic}, we prepare some lemmas which will play 
key roles. In \ref{subsect:wit}, we shall give some conditions 
under which weak index theorem holds.
In particular, we give an effective bound for $m$ 
such that $R^i p_{Y*}(p_X^*(E(mH)) \otimes {\cal E})=0$, $i>0$.
By using these results, we discuss the problem on the preservation 
of the stability conditions.
 
In section \ref{sect:birat}, we consider birational properties
of moduli spaces on abelian surfaces $X$ with $\rho(X)=1$.
We shall show that the Fourier-Mukai transform induced by
the Poincar\'{e} line bundle induces a birational 
correspondence which was conjectured in
\cite{Y:7} (Theorem \ref{thm:birat1}).

Finally we would like to mention that 
Verbitsky \cite{V:1} gets some opposite results to our
results.

\section{Preliminaries}\label{sect:pre}

Let $X$ be a K3 surface or an abelian surface defined over ${\Bbb C}$.
We define a lattice structure $\langle \quad,\quad \rangle$
on 
$H^{ev}(X,{\Bbb Z}):=\bigoplus_{i=0}^2 H^{2i}(X,{\Bbb Z})$
by 
\begin{equation}
\begin{split}
\langle x,y \rangle:=&-\int_X x^{\vee} \wedge y\\
=& \int_X(x_1 \wedge y_1-x_0 \wedge y_2-x_2 \wedge y_0),
\end{split} 
\end{equation}
where $x_i \in H^{2i}(X,{\Bbb Z})$ (resp. $y_i \in H^{2i}(X,{\Bbb Z})$)
is the $2i$-th component of $x$ (resp. $y$)
and $x^{\vee}=x_0-x_1+x_2$.
It is now called the Mukai lattice.
For a coherent sheaf $E$ on $X$,
\begin{equation}
\begin{split}
v(E):=&\ch(E) \sqrt{\td_X}\\
=&\rk(E)+c_1(E)+(\chi(E)-\epsilon \rk(E))\varrho_X \in H^{ev}(X,{\Bbb Z})
\end{split}
\end{equation}
is called the Mukai vector of $E$, where 
$\epsilon=0,1$ according as $X$ is an abelian surface or a K3 surface
and $\varrho_X$ is the fundamental class of $X$.

In \cite{Y:7}, we introduced the notion of twisted stability.
Let $K(X)$ be the Grothendieck group of $X$.
We fix an ample divisor $H$ on $X$.
For $G \in K(X) \otimes{\Bbb Q}$ with $\rk G>0$, we define the 
$G$-twisted rank, 
degree, and Euler characteristic of $x \in K(X) \otimes {\Bbb Q}$ by
\begin{equation}
\begin{split}
\rk_{G}(x)&:=\rk(G^{\vee} \otimes x),\\
\deg_{G}(x)&:=\deg(G^{\vee} \otimes x)=(c_1(G^{\vee} \otimes x),H),\\
\chi_{G}(x)&:=\chi(G^{\vee} \otimes x).
\end{split}
\end{equation}
We define the $G$-twisted stability as follows.
\begin{defn}
Let $E$ be a torsion free sheaf on $X$.
$E$ is $G$-twisted semi-stable (resp. stable) with respect to $H$, if
\begin{equation}
\frac{\chi_G(F(nH))}{\rk_G(F)} \leq \frac{\chi_{G}(E(nH))}{\rk_{G}(E)},
n \gg 0
\end{equation}
for $0 \subsetneq F \subsetneq E$
(resp. the inequality is strict).
\end{defn}

For a Mukai vector $v$,
we denote the moduli stack of $G$-twisted semi-stable sheaves $E$
with $v(E)=v$ by ${\cal M}_H^G(v)^{ss}$
and the open substack consisting of $G$-twisted stable sheaves
by ${\cal M}_H^G(v)^{s}$.
Let $\overline{M}_H^G(v)$ be the moduli space of $S$-equivalence 
classes of $G$-twisted semi-stable sheaves $E$
with $v(E)=v$.
For a coherent sheaf $E$ on $X$, let $0 \subset F_1 \subset F_2 \subset
\dots \subset F_s=E$
be the Harder-Narasimhan filtration of $E$
with respect to the $\mu$-semi-stability.
We set
\begin{equation}
\begin{split}
\mu_{\max,G}(E)&:=\frac{\deg_G(F_1)}{\rk_G(F_1)}=
 \frac{\deg(F_1)}{\rk F_1}-\frac{\deg G}{\rk G},\\
\mu_{\min,G}(E)&:=\frac{\deg_G(F_s/F_{s-1})}{\rk_G(F_s/F_{s-1})}=
\frac{\deg(F_s/F_{s-1})}{\rk (F_s/F_{s-1})}-
\frac{\deg G}{\rk G}.
\end{split}
\end{equation}

\begin{defn}\label{defn:general1}
Let $v$ be a Mukai vector with $\rk v>0$.
A polarization $H$ on $X$ is general with respect to $v$, if
for every $\mu$-semi-stable sheaf $E$ with $v(E)=v$ and a  
subsheaf $F \ne 0$ of $E$, 
\begin{equation}
\frac{(c_1(F),H)}{\rk F}=\frac{(c_1(E),H)}{\rk E}
\text{ if and only if }
\frac{c_1(F)}{\rk F}=\frac{c_1(E)}{\rk E}.
\end{equation}
\end{defn}

Let $v_0:=r_0+\xi_0+a_0 \varrho_X, r_0>0, \xi_0 \in \NS(X)$ be a 
primitive isotropic Mukai vector on $X$.
We take a general ample divisor $H$
with respect to $v_0$.
We set $Y:=M_H(v_0)$. 
Then $Y$ is an abelian surface (resp. a K3 surface),
if $X$ is an abelian surface (resp. a K3 surface).

By the proof of \cite[Lem. 2.1]{Y:8}, the following lemma holds.
\begin{lem}\cite[Lem. 2.1]{Y:9}\label{lem:class-v_1}
Assume that $H$ is general with respect to $v_0$.
\begin{enumerate}
\item
If $Y$ contains a non-locally free sheaf,
then there is an exceptional vector bundle $E_0$ such that
$v_0=\rk(E_0)v(E_0^{\vee})-\varrho_X$.
Moreover $Y \cong X$ and a universal family is given by
\begin{equation}\label{eq:E_0}
{\cal E}:=\ker(E_0^{\vee} \boxtimes E_0 \to {\cal O}_{\Delta}).
\end{equation}
\item
If $Y$ consists of locally free sheaves, then 
they are $\mu$-stable.
\end{enumerate}
\end{lem}

If $X$ is an abelian surface, then $Y$ consists of $\mu$-stable
vector bundles.
Assume that there is a universal family ${\cal E}$ on
$X \times Y$.
Let $p_X:X \times Y \to X$ (resp. $p_Y:X \times Y \to Y$) be the projection.
We define
${\cal F}_{\cal E}:{\bf D}(X) \to {\bf D}(Y)$ by
\begin{equation}
{\cal F}_{\cal E}(x):={\bf R}p_{Y*}({\cal E} \otimes p_X^*(x)),
x \in {\bf D}(X),
\end{equation}
and
$\widehat{\cal F}_{\cal E}:{\bf D}(Y) \to {\bf D}(X)$ by
\begin{equation}
\widehat{\cal F}_{\cal E}(y):={\bf R}\Hom_{p_{X}}({\cal E}, p_Y^*(y)),
y \in {\bf D}(Y),
\end{equation}
where $\Hom_{p_{X}}(-,-)=p_{X*}{\cal H}om_{{\cal O}_{X \times Y}}(-,-)$
is the sheaf of relative homomorphisms. 
Bridgeland \cite{Br:2} showed that ${\cal F}_{\cal E}$ is an equivalence of
categories and the inverse is given by $\widehat{\cal F}_{\cal E}[2]$.
${\cal F}_{\cal E}$ is now called the Fourier-Mukai functor.
We denote the $i$-th cohomology sheaf $H^i({\cal F}_{\cal E}(x))$ by
${\cal F}_{\cal E}^i(x)$.
${\cal F}_{\cal E}$ also induces an isometry of
the Mukai lattices ${\cal F}_{\cal E}:H^{ev}(X,{\Bbb Z})
\to H^{ev}(Y,{\Bbb Z})$.
We are also interested in the composition of ${\cal F}_{\cal E}$ and the
``taking-dual'' functor ${\cal D}_{Y}:{\mathbf D}(Y) \to
{\mathbf D}(Y)_{op}$
sending $x \in {\mathbf D}(Y)$ to 
${\mathbf R}{\cal H}om(x,{\cal O}_{Y})$,
where ${\mathbf D}(Y)_{op}$ is the opposite category of 
${\mathbf D}(Y)$.
By Grothendieck-Serre duality, 
${\cal G}_{\cal E}:=({\cal D}_{Y} \circ {\cal F}_{\cal E})[2]$
is defined by
\begin{equation}
 {\cal G}_{\cal E}(x):={\mathbf R}\Hom_{p_{Y}}
 ({\cal E} \otimes p_{X}^*(x),{\cal O}_{X \times Y}), 
 x \in {\mathbf D}(X). 
\end{equation}
Let 
$\widehat{{\cal G}}_{\cal E}:{\mathbf D}(Y)_{op} \to 
{\mathbf D}({X})$ be the
inverse of ${\cal G}_{\cal E}$:
\begin{equation}
 \widehat{{\cal G}}_{\cal E}(y):={\mathbf R}\Hom_{p_{{X}}}
 ({\cal E} \otimes p_{Y}^*(y),{\cal O}_{X \times Y}), 
 y \in {\mathbf D}(Y). 
\end{equation}

We set $w_0:=v({\cal E}_{|\{x \}\times Y})=
r_0+\widetilde{\xi}_0+\widetilde{a}_0 \varrho_Y$, 
$x \in X$, $\widetilde{\xi}_0 \in \NS(Y)$. 
For a ${\Bbb Q}$-line bundle $L \in K(X) \otimes {\Bbb Q}$,
we define
\begin{equation}\label{eq:det}
\widehat{L}:=\det(p_{Y!}({\cal F}_{\cal E}(E)))
\end{equation}
where $E$ is an element of $K(X) \otimes {\Bbb Q}$
with $\det E=L^{\vee}$ and
$v(E)=-c_1(L)+\frac{1}{r_0}(c_1(L),\xi_0)\varrho_X$.
More precisely, $c_1(\widehat{L}) \in \NS(Y) \otimes {\Bbb Q}$
is well defined.
We usually identify a ${\Bbb Q}$-divisor class $D$ with
the ${\Bbb Q}$-line bundle ${\cal O}_X(D)$.
Hence $\widehat{D}$ denotes a ${\Bbb Q}$-divisor on $Y$
such that ${\cal O}_Y(\widehat{D})=\widehat{{\cal O}_Y(D)}$.
Since $H$ is general with respect to $v_0$,
$\widehat{H}$ is an ample divisor (cf. \cite{Y:7}).
Every Mukai vector $v$ can be uniquely written as 
\begin{equation}
 v=l v_0^{\vee}+a \varrho_X+d(H-\frac{1}{r_0}(H,\xi_0)\varrho_X)+(D
 -\frac{1}{r_0}(D,\xi_0)\varrho_X),
\end{equation}
where $l, a, d \in {\Bbb Q}$,
and $D \in \NS(X)\otimes {\Bbb Q} \cap H^{\perp}$.
It is easy to see that
$l=-\langle v,\varrho_X \rangle/\rk v_0$,
$a=-\langle v, v_0^{\vee} \rangle/\rk v_0$ and 
$d=\deg_{G_1}(v)/(\rk v_0(H^2))$,
where $G_1:={\cal E}_{|X \times \{y \}}^{\vee}$
for a point $y \in Y$.
Since ${\cal F}_{\cal E}(v_1)=\varrho_Y$ and
$\widehat{{\cal F}}_{\cal E}(w_1)=\varrho_X$,
we get
\begin{equation}\label{eq:deg-preserve}
{\cal F}_{\cal E}(l v_0^{\vee}+a \varrho_X+
(dH+D-\frac{1}{r_0}(dH+D,\xi_0)\varrho_X))=
l \varrho_Y+a w_0-
(d\widehat{H}+\widehat{D}+
\frac{1}{r_0}(d\widehat{H}+\widehat{D},\widetilde{\xi}_0)\varrho_Y)
\end{equation}
where $\widehat{D} \in \NS(X) \otimes {\Bbb Q} \cap \widehat{H}^{\perp}$.

Throughout this note, we assume the following two conditions:
\begin{itemize}
\item[(\#1)]
$\widehat{H}$ is general with respect to $w_1$.
\item[(\#2)]
${\cal E}_{|\{x \} \times Y}$ is stable with respect to $\widehat{H}$.
\end{itemize}
\begin{rem}\label{rem:00/08/17}
The assumption $(\#1,2)$ holds for all general $H$, if one of the 
following conditions holds: 
\begin{enumerate}
\item
$X$ is an abelian surface, 
\item
$\NS(X)\cong {\Bbb Z}$,
\item
$Y$ consists of non-locally free sheaves.
\end{enumerate}
For another example, see \cite{BBH:1}.
\end{rem}

Since ${\cal F}_{\cal E}$ is an equivalence of categories,
we get the following.
\begin{lem}\label{lem:spectral}
Let $E$ be a coherent sheaf on $Y$.
Then we have a spectral sequence
\begin{equation}
E_2^{p,q}={\cal F}_{\cal E}^p(\widehat{{\cal F}}_{\cal E}^q(E)) 
\Rightarrow
E_{\infty}^{p+q}=
\begin{cases}
E, &p+q=2,\\
0, &p+q \ne 2.
\end{cases}
\end{equation}
In particular,
\begin{enumerate}
\item
${\cal F}_{\cal E}^p(\widehat{\cal F}_{\cal E}^0(E))=0$, $p=0,1$. 
\item
${\cal F}_{\cal E}^p(\widehat{\cal F}_{\cal E}^2(E))=0$, $p=1,2$. 
\item
There is an injective homomorphism
${\cal F}_{\cal E}^0(\widehat{\cal F}_{\cal E}^1(E)) \to
{\cal F}_{\cal E}^2(\widehat{\cal F}_{\cal E}^0(E))$.
\end{enumerate}
\end{lem}

\section{Counter examples}\label{sect:counter}

In this section, we show that the Fourier-Mukai transform 
does not always preserve the $\mu$-stability of vector bundles
even if $Y$ consists of $\mu$-stable vector bundles.

\subsection{Example 1: an abelian surface case}\label{subsect:counter1}

We shall first give an example for an abelian surface $X$ with
$\rho(X) \geq 2$. We shall treat the $\rho(X)=1$ case later
(see Lemma \ref{lem:counter}). 
Let $(X,H)$ be a polarized abelian surface and $\widehat{X}$ the dual of $X$.
Let ${\cal P}$ be the Poincar\'{e} line bundle on
$X \times \widehat{X}$.
Under the natural identification $H_2(X,{\Bbb Z})=H^2(\widehat{X},{\Bbb Z})$,
$\widehat{D} \in H^2(\widehat{X},{\Bbb Z})$ 
denotes the Poincar\'{e} dual of $D \in H^2(X,{\Bbb Z})$.
This notation is compatible with \eqref{eq:det}.
We shall show that 
there is a Mukai vector $v$ such that
${\cal F}_{\cal P}$ does not preserve
the stability for all $E \in M_H(v)$.
 
Assume that there is a divisor $D$
such that $(D,H)=1$ and $(D^2)=-2$, (and hence we assume that
$\rho(X) \geq 2$). 
Then $M_H((r+1)+D-\varrho_X)$ consists of $\mu$-stable sheaves.
By \cite[Prop. 3.5]{Y:7},
we have an isomorphism
\begin{equation}
 {\cal F}_{\cal P}:M_H((r+1)+D-\varrho_X) \to 
M_{\widehat{H}}(1+\widehat{D}-(r+1)\varrho_{\widehat{X}}). 
\end{equation}
Let $I_Z(\widehat{D})$ be an element of 
$M_{\widehat{H}}(1+\widehat{D}-(r+1)\varrho_{\widehat{X}})
=X \times \Hilb_{\widehat{X}}^r$.
Applying $\widehat{{\cal F}}_{\cal P}$ to the exact sequence
\begin{equation}
 0 \to I_Z(\widehat{D}) \to {\cal O}_{\widehat{X}}(\widehat{D}) \to
 {\cal O}_Z \to 0,
\end{equation}
we get an exact sequence
\begin{equation}
 0 \to P \to E \to F
 \to 0,
\end{equation}
where $P:=p_{X*}({\cal P}^{\vee} \otimes p_{\widehat{X}}^*({\cal O}_Z))$,
$E:=\widehat{{\cal F}}_{\cal P}^1(I_Z(\widehat{D}))$ and
$F:=\widehat{{\cal F}}_{\cal P}^1({\cal O}_{\widehat{X}}(\widehat{D}))$
is a line bundle with 
$v(F)=1+D-\varrho_X$.
By this exact sequence, $E$ is locally free.
We consider the Fourier-Mukai transform of 
a $\mu$-stable vector bundle $E(D)$.
Since $P(D)$ and 
$F(D)$ satisfy $\IT_1$,  
$E(D)$ also satisfies $\IT_1$
and ${\cal F}_{\cal P}^1(E(D))$
fits in an exact sequence
\begin{equation}
0 \to {\cal F}_{\cal P}^1(P(D)) \to 
{\cal F}_{\cal P}^1(E(D)) \to {\cal F}_{\cal P}^1(F(D)) \to 0.
\end{equation}
Since $v( {\cal F}_{\cal P}^1(P(D)) )=r(1+\widehat{D}-\varrho_{\widehat{X}})$ 
and
$v({\cal F}_{\cal P}^1(F(D)))=4+2\widehat{D}-\varrho_{\widehat{X}}$,
we get that
\begin{equation}
\frac{(c_1({\cal F}_{\cal P}^1(P(D))),\widehat{H})}{r}=1>
\frac{1}{2}=\frac{(c_1({\cal F}_{\cal P}^1(F(D))),\widehat{H})}{4}.
\end{equation}
Therefore ${\cal F}_{\cal P}^1(E(D))$ is not $\mu$-semi-stable.

\begin{rem}
Assume that $X$ is a product of two elliptic curves $C_1,C_2$:
$X=C_1 \times C_2$.
We set $f:=\{0 \}\times C_2$ and $g:=C_1 \times \{0 \}$.
Then $H:=2f+g$ and $D:=-f+g$ satisfy the above conditions.
\end{rem}

\subsection{Example 2: a K3 surface case}\label{subsect:counter2}

Let $(X,H)$ be a polarized K3 surface such that
$\Pic(X)={\Bbb Z}H$ with $(H^2)=2n$.
Then $v_0:=k^2n+kH+\varrho_{X}$, $k>0$ is a primitive isotropic Mukai vector.
We assume that $kH$ is very ample.
\begin{lem}
$M_H(v_0) \cong M_H(1+kH+k^2n \varrho_{X}) \cong X$
and $M_H(v_0)$ consists of $\mu$-stable vector bundles.
\end{lem}

\begin{proof}
We use the Fourier-Mukai functor 
${\cal F}_{I_{\Delta}}:{\bf D}(X) \to {\bf D}(X)$, where 
$I_{\Delta}$ is the ideal sheaf of the diagonal 
$\Delta \subset X \times X$.
Since $kH$ is very ample, $\IT_0$ holds for $I_s$, $s \in X$ and  
${\cal F}_{I_{\Delta}}^0(I_s(kH))$ is a simple vector bundle with
the Mukai vector $v_0^{\vee}$. Since $\Pic(X)={\Bbb Z}H$,
it is also stable (\cite[Prop. 3.14]{Mu:4}). 
Moreover ${\cal F}_{I_{\Delta}}^0(I_s(kH))$
is $\mu$-stable: Indeed
let 
\begin{equation}
0 \subset F_1 \subset F_2 \subset \dots \subset 
F_s={\cal F}_{I_{\Delta}}^0(I_s(kH))
\end{equation}
 be the Jordan-H\"{o}lder filtration of 
${\cal F}_{I_{\Delta}}^0(I_s(kH))$ with respect to the $\mu$-stability.
Then we can set that
$v(F_i/F_{i-1})=r_i(kn-H)+a_i \varrho_X$, where
$r_i$ and $a_i$ are integers with $0<r_i \leq k$.
Since $F_i/F_{i-1}$ are $\mu$-stable,
we get $\langle v(F_i/F_{i-1})^2 \rangle=2r_i n(r_i-ka_i) \geq -2$.
If the equality holds, then $n=1$, $r_i=1$ and $k=1,2$.
In these cases, $kH$ is not very ample.
Hence the equality does not hold. Thus $r_i-ka_i \geq 0$.
On the other hand, since $\sum_i r_i=k$ and $\sum_i a_i=1$,
we get that $r_i-ka_i=0$ for all $i$.
Since $r_i \leq k$, we should have $s=1$. Thus 
${\cal F}_{I_{\Delta}}^0(I_s(kH))$ is $\mu$-stable.
 Therefore we get an isomorphism
\begin{equation}
M_H(v_0) \to M_H(v_0^{\vee}) \to M_H(1+kH+k^2n \varrho_{X}) \cong X.
\end{equation}
\end{proof}

\begin{rem}
If $n=1$ and $k=1,2$, then $kH$ is not very ample.
In these cases, we still have isomorphisms $M_H(v_0) \cong X$:
Indeed $\widehat{E}:=
\Ext^2_{\pi_2}(I_{\Delta} \otimes \pi_1^*(I_s(kH)),{\cal O}_X)$, 
$s \in X$ is a stable sheaf with the Mukai vector $v_0$,
where $\pi_i:X \times X \to X$, $i=1,2$ are two projections. 
If $k=1$, then $\widehat{E}$ is isomorphic to 
$I_s(H)$. 
If $k=2$, then $\widehat{E}$ is isomorphic to
$\ker(ev:E_0 \otimes \Hom(E_0,{\Bbb C}_s) \to {\Bbb C}_s)$, $s \in X$, where
$E_0$ is a stable and rigid vector bundle with $v(E_0)=2+H+\varrho_X$.
Therefore it is not $\mu$-stable.
\end{rem}
\vspace{1pc}

Under this identification, we shall construct a universal family
on $X \times X$. 
Let $X_i$, $i=1,2,3$ be three copies of $X$.
Let $p_{ij}: X_1 \times X_2 \times X_3 \to X_i \times X_j$ and
$p_i: X_1 \times X_2 \times X_3 \to X_i$ be the
projections.
We set
\begin{equation}
{\cal E}:=p_{13*}(p_{12}^*(I_{\Delta}) \otimes
p_{23}^*(I_{\Delta}) \otimes p_2^*({\cal O}_{X_2}(kH)))^{\vee}.
\end{equation}
We set $E_0:=q_{3*}(I_{\Delta} \otimes q_1^*({\cal O}_{X_1}(kH)))$,
where $q_i:X_1 \times X_3 \to X_i$, $i=1,3$ are the projections.
Then $E_0$ is a $\mu$-stable vector bundle with 
$\langle v(E_0)^2 \rangle=-2$ and ${\cal E}^{\vee}$ fits in an exact sequence
\begin{equation}\label{eq:exact}
0 \to {\cal E}^{\vee} \to q_3^*(E_0) \to
I_{\Delta} \otimes q_1^*{\cal O}_{X_1}(kH) \to 0.
\end{equation}
Hence ${\cal E}_{|X_1 \times \{s \}}^{\vee}={\cal F}_{I_{\Delta}}^0(I_s(kH))$,
$s \in X_3$ is a $\mu$-stable vector bundle with the Mukai vector $v_0^{\vee}$.
Since ${\cal E}$ is invariant under the natural action of
${\frak S}_2$ on $X_1 \times X_3$,
${\cal E}_{|\{s \} \times X_3}^{\vee}$, $s \in X_1$
 is also a $\mu$-stable vector bundle with
the Mukai vector $v_0^{\vee}$. 
Hence $M_H(v_0) \cong X_1 \cong X_3$ and 
${\cal E}$ becomes a universal family on $X_1 \times X_3$.
By \eqref{eq:exact} and the ${\frak S}_2$-symmetry,
we see that
\begin{equation}\label{eq:F^1=0}
H^1(X_3,{\cal E}_{|\{s\} \times X_3}^{\vee})=
H^1(X_1,{\cal E}_{|X_1 \times \{s\}}^{\vee})=0.
\end{equation}
\begin{rem}
By the exact sequence \eqref{eq:exact}, we see that
${\cal E}_{|\{s\} \times X_3}^{\vee}={\cal F}_{\cal G}^0(I_s)$,
$s \in X_1$,
where 
${\cal G}:=\ker (ev:q_1^*(E_0^{\vee}) \otimes 
q_3^*(E_0) \to {\cal O}_{\Delta})$. 
\end{rem}
\vspace{1pc}

We consider the Fourier-Mukai transform 
${\cal F}_{{\cal E}^{\vee}}:{\bf D}(X_1) \to {\bf D}(X_3)$.
By the construction of ${\cal E}$, we get a decomposition:
\begin{equation}
{\cal F}_{{\cal E}^{\vee}}={\cal F}_{I_{\Delta}(q_2^{-1}(kH))} \circ
{\cal F}_{I_{\Delta}}:{\bf D}(X_1) \to {\bf D}(X_2) \to {\bf D}(X_3).
\end{equation}
Then we see that the induced homomorphism
$H^{ev}(X_1,{\Bbb Z}) \to H^{ev}(X_3,{\Bbb Z})$ is given by  
\begin{equation}\label{eq:F^H}
\begin{split}
{\cal F}_{{\cal E}^{\vee}}(1) &=1,\\
{\cal F}_{{\cal E}^{\vee}}(\xi) &=-k(H,\xi)+\xi,\;\; 
\xi \in H^2(X_1,{\Bbb Z}),\\
{\cal F}_{{\cal E}^{\vee}}(\varrho_{X_1}) &=k^2n-kH+\varrho_{X_3}.
\end{split}
\end{equation}

\begin{lem}\label{lem:02-3}
Let $F$ be a stable sheaf on $X_1$ with $c_1(F)=H$ and
$H^0(X_1,F)=0$. Then
$\Hom({\cal E}_{|X_1 \times \{s\}},F)=0$, $s \in X_3$.
\end{lem}
\begin{proof}
Assume that there is a non-zero map 
$\phi:{\cal E}_{|X_1 \times \{s\}} \to F$.
By the stability of ${\cal E}_{|X_1 \times \{s\}}$,
$c_1(\im \phi)=lH, l>0$.
Since $c_1(F)=H$ and $F$ is stable, we see that
$F/\im \phi$ is of $0$-dimensional.
Thus $\phi$ is surjective in codimension 1.
Since there is an exact sequence
\begin{equation}
{\cal O}_{X_1}^{\oplus (k^2n+1)} \to
{\cal E}_{|X_1 \times \{s\}} \to {\Bbb C}_s \to 0,
\end{equation}
we have a generically surjective map
${\cal O}_{X_1}^{\oplus (k^2n+1)} \to F$.
Hence $H^0(X_1,F) \ne 0$, which is a contradiction.
\end{proof}

\begin{lem}\label{lem:02-1}
For $E \in M_H(r+H-a \varrho_{X_1})$ with $\dim H^0(X_1,E)=t$,
$t \leq r$, we have 
$\Hom({\cal E}_{|X_1 \times \{s\}},E)=0$.
\end{lem}

\begin{proof}
Since $t \leq r$, we have an exact sequence
\begin{equation}
0 \to {\cal O}_{X_1}^{\oplus t} \to E \to F \to 0,
\end{equation}
where $F$ is a stable sheaf with $c_1(F)=H$ (\cite[Lem. 2.1]{Y:5}).
Since $H^0(X_1,F)=0$ and 
$\Hom({\cal E}_{|X_1 \times \{s\}},{\cal O}_{X_1})=0$,
we get our claim.
\end{proof}

\begin{lem}\label{lem:02-2}
$\IT_2$ holds for ${\cal O}_{X_1}$ 
with respect to ${\cal F}_{{\cal E}^{\vee}}$.
\end{lem}

\begin{proof}
Since $c_1({\cal E}_{|X_1 \times \{s\}})=kH$, $s \in X_3$,
we get that 
$\Hom({\cal E}_{|X_1 \times \{s\}},{\cal O}_{X_1})=0$ for $s \in X_3$.
By \eqref{eq:F^1=0}, we get that
$\Ext^1({\cal E}_{|X_1 \times \{s\}},{\cal O}_{X_1})=0$.
Therefore the claim holds.
\end{proof}

\begin{prop}
Assume that $kn>r$.
Then for $E \in M_H(r+H-a \varrho_{X_1})$ with $\dim H^0(X_1,E) \leq r$,
$\IT_1$ holds with respect to ${\cal F}_{{\cal E}^{\vee}}$.
Moreover if $1+k(r+a)<k^2n$ and $H^1(X_1,E) \ne 0$, then 
${\cal F}_{{\cal E}^{\vee}}^1(E)$ is not $\mu$-semi-stable.
\end{prop}

\begin{proof}
By Lemma \ref{lem:02-1}, $\Hom({\cal E}_{|X_1 \times \{s\}},E)=0$
for all $s \in X_3$.
Since $kn>r$, the stability condition implies that 
$\Ext^2({\cal E}_{|X_1 \times \{s\}},E) \cong 
\Hom(E,{\cal E}_{|X_1 \times \{s\}})^{\vee}=0$.
Therefore $\IT_1$ holds.

If $\Ext^1(E,{\cal O}_{X_1})=H^1(X_1,E)^{\vee} \ne 0$, 
then we consider the universal extension
\begin{equation}
0 \to V \otimes {\cal O}_{X_1}
\to E' \to E \to 0,
\end{equation}
where $V=\Ext^1(E,{\cal O}_{X_1})^{\vee}$.
Since $\dim H^0(X_1,E) \leq r$, by using the Riemann-Roch theorem,
we get $\dim V \leq a$,
and hence 
$\rk E' \leq r+a$.
If $k^2 n>(a+r)k+1$, then 
$1/\rk E' \geq 1/(r+a) >1/(kn)$.
Since $E'$ is stable (\cite[Cor. 2.2]{Y:5}),
$\Ext^2({\cal E}_{|X_1 \times \{s \}},E')=0$. 
Therefore $\IT_1$ holds  
for $E'$.
By Lemma \ref{lem:02-2}, $\IT_2$ holds for ${\cal O}_{X_1}$.
Thus we have an exact sequence
\begin{equation}
 0 \to {\cal F}_{{\cal E}^{\vee}}^1(E') \to
 {\cal F}_{{\cal E}^{\vee}}^1(E) \to
 V \otimes {\cal F}_{{\cal E}^{\vee}}^2({\cal O}_{X_1})
 \to 0.
\end{equation}
By \eqref{eq:F^H}, we see that
\begin{equation}
 \begin{split}
  v({\cal F}_{{\cal E}^{\vee}}^2({\cal O}_{X_1}))&=(k^2n+1)-kH+\varrho_{X_3},\\
  v({\cal F}_{{\cal E}^{\vee}}^1(E))&=((ak^2+2k)n-r)-(ak+1)H+a \varrho_{X_3}.
 \end{split}
\end{equation}
Hence 
\begin{equation}
 \mu({\cal F}_{{\cal E}^{\vee}}^1(E))-
 \mu(V \otimes {\cal F}_{{\cal E}^{\vee}}^2({\cal O}_{X_1}))=
 \frac{-2n(ak+1)}{(ak^2+2k)n-r}+\frac{2nk}{k^2n+1}=
 \frac{2n(k^2n-((a+r)k+1))}{((ak^2+2k)n-r)(k^2n+1)}>0.
\end{equation}
Therefore ${\cal F}_{{\cal E}^{\vee}}^1(E)$ is not $\mu$-semi-stable.
\end{proof}

If $a>r$, then $H^0(X_1,E)=0$ for a general $E \in M_H(r+H-a \varrho_{X_1})$
(\cite{Y:5}).
Hence we get the following corollary.
\begin{cor}\label{cor:counter}
If $kn \geq a >r$, then for a general element $E \in M_H(r+H-a \varrho_{X_1})$,
$\IT_1$ holds with respect to ${\cal F}_{{\cal E}^{\vee}}$.
Moreover if $1+k(r+a)<k^2n$, then 
${\cal F}_{{\cal E}^{\vee}}^1(E)$ is not $\mu$-semi-stable.
\end{cor}

\begin{rem}
In the above example, we used the moduli space of
$\mu$-stable vector bundles.
For the Fourier-Mukai transform induced by a moduli space
consisting of non-locally free sheaves, it is much easy to
construct such an example:   
For the same $E$ in Corollary \ref{cor:counter},
$\WIT_1$ holds with respect to ${\cal F}_{I_{\Delta}}$ and
${\cal F}_{I_{\Delta}}^1(E)$ fits in an exact sequence
\begin{equation}
 0 \to E \to {\cal F}_{I_{\Delta}}^1(E) \to H^1(X,E) \otimes {\cal O}_X
 \to 0. 
\end{equation}
Hence ${\cal F}_{I_{\Delta}}$ does not preserve the stability condition.
\end{rem}

\section{Asymptotic results}\label{sect:asymptotic}
We keep the notation in section \ref{subsect:counter1}.
For a sem-stable sheaf $E$ on an abelian surface 
$X$ with $v(E)=r+\xi+a \varrho_X$, $\xi \in \NS(X)$
and a subsheaf $E_1$ with $v(E_1)=r_1+\xi_1+a_1 \varrho_X$,
$\xi_1 \in \NS(X)$,
we see that
\begin{equation}
\begin{split}
&\frac{\deg({\cal F}_{\cal P}(E_1(mH)))}{\rk({\cal F}_{\cal P}(E_1(mH)))}
-\frac{\deg({\cal F}_{\cal P}(E(mH)))}{\rk({\cal F}_{\cal P}(E(mH)))} \\
=&
\frac{-(\xi_1+mr_1 H,H)}{\chi(E_1(mH))}-
\frac{-(\xi+mr H,H)}{\chi(E(mH))} \\
=&\frac{(r \xi_1-r_1 \xi,H)m^2(H^2)/2+
(ra_1-r_1a)m(H^2)+((\xi,H)a_1-(\xi_1,H)a)}{\chi(E_1(mH))\chi(E(mH))}
\end{split}
\end{equation}
and
\begin{equation}
\frac{\chi({\cal F}_{\cal P}(E_1(mH)))}{\rk({\cal F}_{\cal P}(E_1(mH)))}
-\frac{\chi({\cal F}_{\cal P}(E(mH)))}{\rk({\cal F}_{\cal P}(E(mH)))}
=\frac{r_1\chi(E(mH))-r\chi(E_1(mH))}{\chi(E_1(mH))\chi(E(mH))}.
\end{equation}
Hence if $m$ is sufficiently large, then $E_1$ does not induce a 
destabilizing subsheaf of ${\cal F}_{\cal P}(E(mH))$. 

In this section, we consider 
the preservation of stability
for $E(mH)$, $m \gg 0$ under the Fourier-Mukai transform 
${\cal F}_{\cal E}:{\bf D}(X) \to {\bf D}(Y)$
for a general ${\cal E}$ with the conditions $(\#1,2)$.

\subsection{Basic lemmas}\label{subsect:basic}
Keep the notation in section \ref{sect:pre}.
We set $(H^2)=2n$.
We assume that ${\cal E}$ satisfies conditions $(\#1,2)$.  
We set 
$G_1:={\cal E}_{|X \times \{y \}}^{\vee}$ and
$G_2:={\cal E}_{|\{x \} \times Y}$ for some $x \in X$ and $y \in Y$.

We note that 
\begin{equation}\label{eq:deg}
d:=\frac{\deg_{G_1}(E)}{r_0 (H^2)}
 \in \frac{1}{2nr_0}{\Bbb Z}.
\end{equation}

\begin{lem}\label{lem:slope-bound}
Let $E$ be a $\mu$-semi-stable sheaf with $\deg_{G_1}(E)>0$.
\begin{enumerate}
\item[(1)]
Assume that $E$ satisfies $\IT_0$ with respect to ${\cal F}_{\cal E}$.
Then 
\begin{equation}
\max \{\deg_{G_2}(F)|\;F \subset {\cal F}_{\cal E}^0(E) \}<0.
\end{equation}
\item[(2)]
Assume that $E$ satisfies $\WIT_2$ with respect to ${\cal G}_{\cal E}$.
Then 
\begin{equation}
\min \{\deg_{G_2^{\vee}}(G)| \;{\cal G}_{\cal E}^2(E) \to G \to 0 \}>0.
\end{equation}
\end{enumerate}
\end{lem}

\begin{proof}
We shall only prove (1). The proof of (2) is similar. 
If the claim does not hold, then
there is an exact sequence
\begin{equation}
0 \to F_1 \to {\cal F}_{\cal E}^0(E) \to F_2 \to 0
\end{equation}
such that $F_1$ is a torsion free sheaf with 
$\mu_{\min,G_2}(F_1) \geq 0$ and
$F_2$ is a torsion free sheaf with $\mu_{\max,G_2}(F_2) < 0$.
Applying $\widehat{\cal F}_{\cal E}$ to this exact sequence, 
we get a long exact sequence
\begin{equation}\label{eq:FM1}
\begin{CD}
0 @>>> \widehat{\cal F}_{\cal E}^0(F_1) @>>> 0 @>>> 
\widehat{\cal F}_{\cal E}^0(F_2)\\
@>>>\widehat{\cal F}_{\cal E}^1(F_1) @>>> 0 @>>> 
\widehat{\cal F}_{\cal E}^1(F_2)\\
@>>> \widehat{\cal F}_{\cal E}^2(F_1) @>>> E @>>> 
\widehat{\cal F}_{\cal E}^2(F_2) @>>> 0
\end{CD}
\end{equation}
Since $\mu_{\max,G_2}(F_2) < 0$, we have $\widehat{\cal F}_{\cal E}^0(F_2)=0$.
By $\mu_{\min,G_2}(F_1) \geq 0$, we see that 
$\widehat{\cal F}_{\cal E}^2(F_1)$ is $0$-dimensional or $0$.
Since $E$ is torsion free or purely 1-dimensional,
$\widehat{\cal F}_{\cal E}^2(F_1) \to E$ is a 0-map.
Hence $\widehat{\cal F}_{\cal E}^1(F_2) \cong 
\widehat{\cal F}_{\cal E}^2(F_1)$
satisfies $\IT_0$.
By Lemma \ref{lem:spectral},
we have an injection ${\cal F}_{\cal E}^0(\widehat{\cal F}_{\cal E}^1(F_2))
\to {\cal F}_{\cal E}^2(\widehat{\cal F}_{\cal E}^0(F_2))=0$.
Hence $\widehat{\cal F}_{\cal E}^1(F_2)=0$, which implies that
$F_1=0$. Therefore our claim holds.
\end{proof}

\begin{lem}\label{lem:sing}
Under the same assumptions, 
$\Hom({\cal E}_{|\{x \} \times Y },
{\cal F}_{\cal E}^0(E))=0$ for all $x \in X$
and the set
\begin{equation}
S:=\{x \in X|\Ext^1({\cal E}_{|\{x \} \times Y},
{\cal F}_{\cal E}^0(E)) \ne 0 \}
\end{equation}
coincides with the set 
\begin{equation}
\{x \in X|\text{ $E \otimes {\cal O}_{X,x}$ is not free } \},
\end{equation}
where $ {\cal O}_{X,x}$ is the stalk of ${\cal O}_X$ at $x$.
\end{lem}

\begin{proof}
By Lemma \ref{lem:slope-bound}, $\mu_{\max,G_2}({\cal F}_{\cal E}^0(E))<0$.
Hence the first claim holds.
Since $\WIT_2$ holds for ${\cal F}_{\cal E}^0(E)$,
$S$ is a proper subset of $X$.
By the base change theorem, we get our claim.
\end{proof} 

\begin{lem}\label{lem:torsion-free}
If $\Hom({\cal E}_{|\{x \} \times Y},F)=0$ 
for all $x \in X$, then
$\widehat{{\cal F}}_{\cal E}^1(F)=\Ext^1_{p_X}({\cal E},p_{Y}^*(F))$
is locally free.
\end{lem}

\begin{proof}
Since $\Hom({\cal E}_{|\{x \} \times Y},F)=0$ 
for all $x \in X$, 
there is a complex of vector bundles
$V_1 \to V_2$ such that
$\widehat{{\cal F}}_{\cal E}^i(F)$, $i=1,2$ are
cohomology sheaves.
Since $X$ is a smooth surface,
$\widehat{{\cal F}}_{\cal E}^1(F)=
\Ext^1_{p_X}({\cal E},p_{Y}^*(F))$
is locally free.   
\end{proof}
The following lemma and its variants will play important roles 
in subsections \ref{subsect:wit} and \ref{subsect:asymptotic}.

\begin{lem}\label{lem:key}
We set $v:=lv_0^{\vee}+a \varrho_X+
(dH+D)-(dH+D,\xi_0)\varrho_X/r_0 \in H^{ev}(X,{\Bbb Z})$,
where $l,a >0$ and
$D \in \NS(X) \otimes {\Bbb Q} \cap H^{\perp}$.
We set $N:=\max\{4r_0^3 l^2+1/(2n),2r_0^2 l(\langle v^2 \rangle-(D^2)) \}$.
Then the following hold: 
\begin{enumerate}
\item[(1)]
If $d > N$, then for any
$G_2$-twisted stable sheaf $F_1$ with
\begin{equation}
\text{$v(F_1)=a_1 w_0+l_1 \varrho_Y-(d_1 \widehat{H}+\widehat{D}_1
+
%\frac{1}{r_0}
(d_1 \widehat{H}+\widehat{D}_1,\widetilde{\xi}_0)\varrho_Y/r_0)$,
$0<d_1<d$ and
$d_1/a_1 \leq d/a$}, 
\end{equation}
we have $l_1 \leq ld_1/d$.
\item[(2)]
If $d > N$, then for any
$G_1$-twisted stable sheaf $E_1$ with
\begin{equation}
\text{
$v(E_1)=l_1 v_0^{\vee}+a_1 \varrho_X+(d_1 {H}+{D}_1
-(d_1 {H}+{D}_1,{\xi}_0)\varrho_X/r_0)$, 
$0<d_1<d$ and
$d_1/l_1 < d/l$,}
\end{equation} 
we have $a_1 < ad_1/d$.
\end{enumerate}
\end{lem}

\begin{proof}
We set $s:=\langle v^2 \rangle/2=-r_0la+d^2n+(D^2)/2$.
We shall first prove (1).
Let $F_1$ be a $G_2$-twisted stable sheaf with
$v(F_1)=a_1 w_0+l_1 \varrho_Y-(d_1 \widehat{H}+\widehat{D}_1
+(d_1 \widehat{H}+\widehat{D}_1,\widetilde{\xi}_0)\varrho_Y/r_0)$,
$0<d_1<d$ and $d_1/a_1 \leq d/a$.
By \eqref{eq:deg}, we get that
$1/(2nr_0) \leq d_1 \leq d-1/(2nr_0)$.
We note that
\begin{equation}\label{eq:F_1}
\begin{split}
\langle v(F_1)^2 \rangle &=2nd_1^2-2l_1 a_1 r_0+(D_1^2)\\
& \leq 2nd_1^2-2l_1r_0 d_1 a/d\\
&=2nd_1^2-2l_1r_0 \frac{d_1}{d}\frac{d^2n-s+(D^2)/2}{r_0 l}\\
&=2nd_1^2-2l_1 d_1 \frac{d^2n-s+(D^2)/2}{dl}.
\end{split}
\end{equation}
We first show that $l_1 < l$ for $d>N$.
Assume that $l_1 \geq l$.
By \eqref{eq:F_1}, we see that
\begin{equation}\label{eq:bd-1}
\begin{split}
-2 \epsilon \leq &\langle v(F_1)^2 \rangle\\
 \leq & 2 d_1^2 n-2 (d^2n-s+(D^2)/2)d_1/d\\
=& 2n d_1\left(d_1-d+\frac{s-(D^2)/2}{dn}\right).
\end{split}
\end{equation}
We set $n_1:=\max\{4r_0+1/(2n r_0), 2r_0(\langle v^2 \rangle-(D^2)) \}$.
We note that $N>n_1$.
We shall show that 
\begin{equation}\label{eq:bd0}
2n d_1\left(d_1-d+\frac{s-(D^2)/2}{dn}\right)
<-2 \epsilon
\end{equation}
for $d > n_1$.
Then by \eqref{eq:bd-1}, we get a contradiction. 
Therefore we have $l_1 < l $ for $d > n_1$. 
\newline
Proof of \eqref{eq:bd0}: 
It is easy to see that \eqref{eq:bd0} follows from the following 
inequality:
\begin{equation}\label{eq:bd1}
d-\frac{s-(D^2)/2}{dn}>
\max\left\{d_1+\frac{\epsilon}{nd_1}\left|
\;d_1=\frac{1}{2nr_0},d-\frac{1}{2nr_0} \right.\right\}
\end{equation}
for all $d \geq n_1$.
Hence we shall show \eqref{eq:bd1}: 
For $d > n_1$,
we have $n(d-1/(2nr_0))>4nr_0$ and
$(s-(D^2)/2)/(dn) < 1/(4nr_0)$.
Hence 
\begin{equation}
d-\frac{1}{2nr_0}+\frac{1}{n(d-1/(2nr_0))} <
d-\frac{1}{2nr_0}+\frac{1}{4nr_0}=d-\frac{1}{4nr_0}
< d-\frac{s-(D^2)/2}{dn}.
\end{equation}
We also get that
$1/(2nr_0)+2r_0 \leq -1/(4nr_0)+1+2r_0<-\frac{s-(D^2)/2}{dn}+d$.
Therefore \eqref{eq:bd1} holds. 

We next show that $l_1 \leq ld_1/d$.
By \eqref{eq:F_1}, we get that
\begin{equation}
-2\epsilon \leq \langle v(F_1)^2 \rangle 
\leq 2nd_1 \left(\left(d_1-\frac{l_1}{l}d \right)+
\frac{l_1}{dnl}(s-(D^2)/2) \right).
\end{equation}
We note that 
\begin{equation}
 2nd_1 \left(\left(d_1-\frac{l_1}{l}d \right)+
 \frac{l_1}{dnl}(s-(D^2)/2) \right)
 <-2 \epsilon
\end{equation}
if and only if
$(l_1/(dnl))(s-(D^2)/2)<
l_1d/l-(d_1+\epsilon/(nd_1))$.
We shall show that $(l_1/l)d-d_1 \leq 0$, if
$d$ is sufficiently large.
Assume that $(l_1/l)d-d_1>0$.
Since $1/(2nr_0) \leq d_1  \leq (l_1/l)d-1/(2nlr_0^2)$,
we get that
\begin{equation}
l_1d/l-(d_1+\epsilon/(nd_1)) \geq
\min\left\{l_1d/l-1/(2nr_0)-\epsilon 2r_0,
1/(2nlr_0^2)-\frac{\epsilon}{n(l_1d/l-1/(2nlr_0)^2)}
\right\}.
\end{equation}
We set $n_2:=4l^2 r_0^3+1/(2n)$. Then
we see that $n(l_1d/l-1/(2nlr_0)) > 4nlr_0^2$ for $d > n_2$.
We set $n_3:=2r_0^2 l+l/(2n)+1/(4n r_0)$. Then we get that
$l_1d/l-1/(2nr_0)-\epsilon 2r_0 \geq 1/(4nlr_0^2)$.
Hence
$l_1d/l-(d_1+\epsilon/(nd_1)) \geq 1/(4nlr_0^2)$ for $d \geq \max\{n_2,n_3\}$.
So if $d>\max\{n_1,n_2,n_3,4r_0^2 l(s-(D^2)/2) \}=N$, then
$\langle v(F_1)^2 \rangle<-2 \epsilon$, which is a contradiction.
Therefore $(l_1/l)d-d_1 \leq 0$ for $d>N$.

We next prove (2).
Assume that $a_1 \geq d_1a/d$.
Since $\langle v(E_1)^2 \rangle=\langle v(F_1)^2 \rangle$,
by the same argument, we get a contradiction.
Therefore $a_1 < d_1a/d$ for $d>N$.
\end{proof}
If $d_1=d$, then we can show the following:
\begin{lem}\label{lem:key-rem}
The same claims in Lemma \ref{lem:key} hold, if $d_1=d$.  That is,
\begin{enumerate}
\item[(1)]
if $d>N$, then 
for a $G_2$-twisted stable sheaf $F_1$ with
\begin{equation}
\text{
$v(F_1)=a_1 w_0+l_1 \varrho_Y-(d_1 \widehat{H}+\widehat{D}_1
+(d_1 \widehat{H}+\widehat{D}_1,\widetilde{\xi}_0)\varrho_Y /r_0)$, 
$d_1=d$ and
$0<d_1/a_1 \leq d/a$,}
\end{equation} 
we have $l_1 \leq l$.
\item[(2)]
If $d > N$, then for any
$G_1$-twisted stable sheaf $E_1$ with
\begin{equation}
\text{$v(E_1)=l_1 v_0^{\vee}+a_1 \varrho_X+(d_1 {H}+{D}_1
-(d_1 {H}+{D}_1,{\xi}_0)\varrho_X /r_0)$, 
$d_1=d$ and
$l_1 > l$, }
\end{equation}
we have $a_1 < a$.
\end{enumerate}
\end{lem}

\begin{proof}
We shall only prove (1).
If $l_1 \geq l+1/r_0$, then we get
\begin{equation}
\begin{split}
\langle v(F_1)^2 \rangle \leq & \langle v^2 \rangle-(D^2)-2a\\
=&\frac{-(2nd^2)+(lr_0+1)(\langle v^2 \rangle-(D^2)) }{(l r_0)}.
\end{split}
\end{equation}
Since $d>4r_0^3 l^2$, we get $nd^2>4n r_0^3 l^2 d$.
Then we see that $nd^2 >4n r_0^3 l^2 d>
(lr_0+1)(\langle v^2 \rangle-(D^2))$ and
$nd^2 >n (4r_0^3 l^2)^2>2lr_0$, and hence $\langle v(F_1)^2 \rangle<-2$.
Therefore we get our claim.
\end{proof}

\begin{lem}\label{lem:key0}
We set $v:=a \varrho_X+
(dH+D)-\frac{1}{r_0}(dH+D,\xi_0)\varrho_X \in H^{ev}(X,{\Bbb Z})$,
where $D \in \NS(X) \otimes {\Bbb Q} \cap H^{\perp}$.
We set $N:=\max \{(\langle v^2 \rangle-(D^2))/2,2r_0+1 \}$.
Then the following hold:
\begin{enumerate}
\item[(1)]
If $a > N$, then 
for any $G_2$-twisted stable sheaf $F_1$ with
\begin{equation}
\text{$v(F_1)=a_1 w_0+l_1 \varrho_Y-(d_1 \widehat{H}+\widehat{D}_1
+(d_1 \widehat{H}+\widehat{D}_1,\widetilde{\xi}_0)\varrho_Y /r_0)$, 
$d_1<d$ and
$0<d_1/a_1 \leq d/a$,}
\end{equation}
 we have
$l_1 \leq 0$.
\item[(2)]
If $a > N$, then 
for any $G_1$-twisted stable sheaf $E_1$ with
\begin{equation}
\text{$v(E_1)=l_1 v_0^{\vee}+a_1 \varrho_X+(d_1 \widehat{H}+\widehat{D}_1
-(d_1 {H}+{D}_1,{\xi}_0)\varrho_X /r_0)$, 
$0<d_1<d$ and $l_1>0$,}
\end{equation}
we have $a_1/d_1 < a/d$.
\end{enumerate}
\end{lem}

\begin{proof}
We shall only prove (1).
Let $F_1$ be a $G_2$-twisted stable sheaf with
\begin{equation}
\text{
$v(F_1)=a_1 w_0+l_1 \varrho_Y-
((d_1 \widehat{H}+\widehat{D}_1)+
(d_1\widehat{H}+\widehat{D}_1,\widetilde{\xi}_0)\varrho_Y /r_0)$, 
$d_1<d$ and $0<d_1/a_1 \leq d/a$}.
\end{equation}
Assume that $l_1>0$. Then $r_0l_1 \geq 1$, and hence we see that
$-2\epsilon \leq \langle v(F_1)^2 \rangle \leq 
d_1^2(H^2)-2r_0l_1a_1 \leq d_1^2(H^2)-2a_1
\leq d_1^2(H^2)-2ad_1/d$.
We set
\begin{equation}
\begin{split}
n_1:= &\max
\left\{d\left(\frac{(H^2)}{2}d_1+\frac{\epsilon}{d_1}\right)\left|
\,\frac{1}{r_0(H^2)} \leq d_1 \leq d-\frac{1}{r_0(H^2)} \right. \right\}\\
=&\max
\left\{d\left(\frac{(H^2)}{2}d_1+\frac{\epsilon}{d_1}\right)\left|
\,d_1=\frac{1}{r_0(H^2)}, d-\frac{1}{r_0(H^2)} \right. \right\}.
\end{split}
\end{equation}
Then we have $d_1^2(H^2)-2ad_1/d<-2 \epsilon$ for $a>n_1$.
Therefore $l_1 \leq 0$ for $a>n_1$.
It is easy to see that
$N:=\max \{(\langle v^2 \rangle-(D^2))/2,2r_0+1 \}>n_1$.
Hence (1) holds. 
\end{proof}

We can also show the following.
\begin{lem}\label{lem:key0-rem}
Keep the notations in Lemma \ref{lem:key0}
\begin{enumerate}
\item[(1)]
If $a > N+1$, then 
for any $G_2$-twisted stable sheaf $F_1$ with
\begin{equation}
\text{$v(F_1)=a_1 w_0+l_1 \varrho_Y-(d_1 \widehat{H}+\widehat{D}_1
+(d_1 \widehat{H}+\widehat{D}_1,\widetilde{\xi}_0)\varrho_Y /r_0)$, 
$d_1=d$ and
$0<d_1/a_1 \leq d/a$,}
\end{equation}
 we have
$l_1 \leq 0$.
\item[(2)]
If $a > N+1$, then 
for any $G_1$-twisted stable sheaf $E_1$ with
\begin{equation}
\text{$v(E_1)=l_1 v_0^{\vee}+a_1 \varrho_X+(d_1 \widehat{H}+\widehat{D}_1
-(d_1 {H}+{D}_1,{\xi}_0)\varrho_X /r_0)$, 
$d_1=d$ and $l_1 >0$,}
\end{equation}
we have $a_1/d_1 < a/d$.
\end{enumerate}
\end{lem}

\begin{cor}\label{cor:slope-stable}
Under the same assumptions as in
Lemma \ref{lem:key0}, or \ref{lem:key},
let $F$ be a $\mu$-semi-stable sheaf
with $v(F)={\cal F}_{\cal E}(v)=a w_0+l \rho_Y
-(d \widehat{H}+\widehat{D}+
(d\widehat{H}+\widehat{D},\widetilde{\xi}_0)\rho_Y /r_0)$. 
Then $F$ is $G_2$-twisted semi-stable.
Moreover, if $F$ is $G_2$-twisted stable, then
it is $\mu$-stable.
\end{cor}

\begin{proof}
Assume that $F$ is not $\mu$-stable.
Let 
\begin{equation}
0 \subset F_1 \subset F_2 \subset \dots \subset F_s=F
\end{equation}
be the Jordan-H\"{o}lder filtration of $F$ with respect to
the $\mu$-stability.
We set $v(F_i/F_{i-1})=a_i w_0+l_i \varrho_Y-(d_i \widehat{H}+\widehat{D}_i
+(d_i \widehat{H}+\widehat{D}_i,\widetilde{\xi}_0)\varrho_Y /r_0)$.
Applying Lemma \ref{lem:key0}, or \ref{lem:key} to each $F_i/F_{i-1}$,
we get that $l_i \leq ld_i/d$.
Then we see that $\sum_i l_i \leq \sum_i ld_i/d =l$.
Since $\sum_i l_i=l$, we have $l_i= ld_i/d$ for all $i$.
Since $d_i/a_i=d/a$, we get $l_i/a_i=l/a$, which implies that
$F$ is $G_2$-twisted semi-stable.
By the same proof, we also see that $F$ is 
$\mu$-stable, provided that $F$ is $G_2$-twisted stable.
\end{proof}

\begin{rem}\label{rem:free}
Under the conditions as in
Lemma \ref{lem:key0-rem}, or \ref{lem:key-rem},
let $F$ be a $\mu$-semi-stable sheaf
with $v(F)={\cal F}_{\cal E}(v)=a w_0+l \rho_Y
-(d \widehat{H}+\widehat{D}+
(d\widehat{H}+\widehat{D},\widetilde{\xi}_0)\rho_Y /r_0)$. 
Then we can easily show that $F$ is locally free.
\end{rem}

\begin{rem}\label{rem:general}
Assume that $l>0$.
We set $w:=a w_0+l \varrho_Y-
(\widehat{\xi}+(\widehat{\xi},\widetilde{\xi}_0)\varrho_Y /r_0)$.
If $(\xi,H)/(H^2) >N$ and
\begin{equation}\label{eq:general}
\min\{-(D^2)|(D,H)=0, D \in \NS(X)\setminus \{0\} \}
>(r_0l)^2(\langle w^2 \rangle+
2(r_0l)^2 \epsilon)/4,
\end{equation}
then $\widehat{H}$ is a general polarization with respect to $w$.

Proof of the claim:
Assume that there is an exact sequence
\begin{equation}
0 \to F_1 \to {\cal F}_{\cal E}^0(E) \to F_2 \to 0
\end{equation}
such that $F_1 (\ne 0)$ is a $\mu$-semi-stable sheaf with
\begin{equation}
v(F_1)=v_1:=a_1 w_0+l_1 \varrho_Y-(\widehat{\xi}_1
+(\widehat{\xi}_1,\widetilde{\xi}_0)\varrho_Y /r_0),
(\xi_1,H)/a_1=(\xi,H)/a
\end{equation}
and $F_2 (\ne 0)$ is a $\mu$-semi-stable sheaf with 
\begin{equation}
v(F_2)=v_2:=a_2 w_0+l_2 \varrho_Y-(\widehat{\xi}_2
+(\widehat{\xi}_2,\widetilde{\xi}_0)\varrho_Y /r_0),
(\xi_2,H)/a_2=(\xi,H)/a.
\end{equation}
By Lemma \ref{lem:key}, we see that $l_1/a_1=l_2/a_2=l/a$,
and $F_1$ and $F_2$ are $G_2$-twisted semi-stable sheaves.
Then we see that $\langle v_i^2 \rangle \geq -2l_i^2r_0^2 \epsilon$.
By a simple calculation, we have an inequality
\begin{equation}
r_0^2l_1l_2(\langle w^2 \rangle+2r_0^2l^2 \epsilon) \geq
-((r_0l_2\xi_1-r_0l_1 \xi_2)^2).
\end{equation}
Since $r_0l_2\xi_1-r_0l_1 \xi_2=(r_0l_2)c_1(\widehat{\cal F}_{\cal E}(v_1))-
(r_0l_1)c_1(\widehat{\cal F}_{\cal E}(v_2)) \in \NS(X)$, we get our claim.
\end{rem}

\subsection{Weak index theorem}
\label{subsect:wit}

We shall give some conditions under which $\WIT_i$ holds
with respect to ${\cal F}_{\cal E}$ or
${\cal G}_{\cal E}$.

\begin{prop}\label{prop:WIT_2}
We set $w:=a w_0+l \varrho_Y-
(d\widehat{H}+\widehat{D}+
(d\widehat{H}+\widehat{D},\widetilde{\xi}_0)\varrho_Y /r_0)$.
Let $F$ be a $G_2$-twisted stable sheaf with 
$v(F)=w$.
If 
\begin{equation}\label{eq:WIT_2}
\frac{d}{r_0 l}>\max \{4lr_0^2+1/(2n r_0 l),
2r_0(\langle w^2 \rangle-(D^2))\},
\end{equation}
then $\WIT_2$ holds for $F$ with respect to $\widehat{{\cal F}}_{\cal E}$
and $\widehat{{\cal F}}_{\cal E}^2(F)$ is torsion free.
\end{prop}

\begin{proof}
By Corollary \ref{cor:slope-stable} and Remark \ref{rem:free}, 
$F$ is a $\mu$-stable vector bundle.
Assume that $\Ext^1({\cal E}_{|\{x_i \} \times Y},F) \ne 0$
for $x_1,x_2,\dots,x_n \in X$.
We take non-zero elements 
$\phi_i \in \Ext^1({\cal E}_{|\{x_i \} \times Y},F)$,
$1 \leq i \leq n$,
and we consider an extension
\begin{equation}
0 \to F \to I \to 
\bigoplus_{i=1}^n {\cal E}_{|\{x_i \} \times Y} \to 0
\end{equation}
whose extension class is given by
$(\phi_1,\phi_2,\dots,\phi_n) \in 
\bigoplus_{i=1}^n \Ext^1({\cal E}_{|\{x_i \} \times Y},F) \cong
\Ext^1(\bigoplus_{i=1}^n{\cal E}_{|\{x_i \} \times Y},F)$.
Let 
\begin{equation}
0 \subset F_1(I) \subset F_2(I) \subset
\dots \subset F_s(I)=I
\end{equation}
 be the Harder-Narasimhan filtration
of $I$ with respect to the $G_2$-twisted semi-stability
(if $s=1$, then $I$ is $G_2$-twisted semi-stable).
We set $I_i:=F_i(I)/F_{i-1}(I)$ and
$v(I_i):=a_i w_0+l_i \varrho_Y-
(d_i\widehat{H}+\widehat{D}_i+
(d_i\widehat{H}+\widehat{D}_i,\widetilde{\xi}_0)\varrho_Y /r_0)$.
Then $-d_1/a_1 \geq -d_2/a_2 \geq \dots \geq -d_s/a_s$.
If $-d_1/a_1 \geq 0$, then $d_1=0$ and 
the natural map $I_1 \to \bigoplus_{i=1}^n {\cal E}_{|\{x_i \} \times Y}$
is injective.
We first assume that ${\cal E}_{|\{x_i \} \times Y}$ is locally free.
Then $I_1$ is also locally free and $I_1$ contains a $\mu$-stable
locally free sheaf $I_1'$ with the same slope as that of $I_1$.
Then $I_1' \cong {\cal E}_{|\{x_i \} \times Y}$ for some $i$,
which is a contradiction.
Therefore we get $d_1 >0$, which also implies that $d_i>0$
for all $i$.
Since $\sum_i d_i=d$, we have $d_i<d$ for all $i$.
Let $I_s'$ be a $\mu$-stable quotient of $I_s$
with the same slope as that of $I_s$.
If $d_s/a_s \geq d/a$, then $d_s/a_s = d/a$ and
we have an injective homomorphism $F \to I_s'$,
which is a contradiction.
Thus $d_s/a_s < d/a$, which implies that $d_i/a_i<d/a$ for all $i$.
By our assumption \eqref{eq:WIT_2}, Lemma \ref{lem:key} (1) implies that
$l_i/l \leq d_i/d$.
Then $\sum_i l_i/l \leq \sum_i d_i/d$.
Since $\sum_i l_i/l=\sum_i d_i/d=1$, we get that $l_i/l=d_i/d$ for all $i$.
Then we see that
\begin{equation}\label{eq:JHF}
\begin{split}
 \sum_i \frac{\langle v(I_i)^2 \rangle-(D_i^2)}{l_i}
 &=\sum_i \frac{2nd_i^2-2r_0l_ia_i}{l_i}\\
 &=\sum_i (2n\frac{d_i}{l_i}d_i-2r_0a_i)\\
 &=2n\frac{d}{l}\sum_i d_i-2r_0 \sum_i a_i\\
 &=\frac{2nd^2-2r_0la}{l}\\
 &=\frac{\langle v(I)^2 \rangle-(D^2)}{l}.
\end{split}
\end{equation}
Since $I_i$ are $G_2$-twisted semi-stable, 
Lemma \ref{lem:bogomolov} below implies that
$\langle v(I_i)^2 \rangle \geq 
-2(r_0 l_i)^2 \epsilon$.
On the other hand, we get $\langle v(I)^2 \rangle-(D^2)=
\langle w^2 \rangle-(D^2)-2nlr_0$.
Hence $n$ is bounded above, which implies that
$\WIT_2$ holds and $\widehat{{\cal F}}_{\cal E}^2(F)$ is torsion free.

We next assume that ${\cal E}_{|\{x_i \} \times Y}$ is not locally free.
Then ${\cal E}_{|\{x_i \} \times Y}=\ker(E_0 \otimes \Hom(E_0,{\Bbb C}_{x_i})
\to {\Bbb C}_{x_i})$ (see \eqref{eq:E_0}). We set $t_0:=\rk E_0$.
We shall show that 
\begin{equation}\label{eq:Ext^1(E_0,F)=0}
 \Ext^1(E_0,F)=0.
\end{equation}
Assume that $\Ext^1(E_0,F) \ne 0$.
We consider a non-trivial extension
\begin{equation}
 0 \to F \to N \to E_0 \to 0.
\end{equation}
Let 
\begin{equation}
0 \subset F_1(H) \subset F_2(N) \subset \dots \subset
F_t(N)=N
\end{equation}
be the Harder-Narasimhan filtration of $N$ with respect to the
$G_2$-twisted semi-stability.
We set 
\begin{equation}
v(F_i(N)/F_{i-1}(N))=a_i w_0+l_i \varrho_Y-
(d_i\widehat{H}+\widehat{D}_i+
(d_i\widehat{H}+\widehat{D}_i,\widetilde{\xi}_0)\varrho_Y /r_0).
\end{equation}
Then $-d_1/a_1 \geq -d_2/a_2 \geq  \dots \geq -d_t/a_t$.
By the same argument as above, we see that
$l_i/l \leq d_i/d$.
Since $\sum_i l_i=l+1/t_0$ and $\sum_i d_i=d$, we get a contradiction.
Therefore $\Ext^1(E_0,F)=0$.
Assume that $I$ is not locally free. 
We set $J:=\{i| \text{ $I$ is not locally free at $x_i$ } \}$
and $K:=\{1,2,\dots, n \} \setminus J$. 
Then we have an exact sequence
\begin{equation}
0 \to F' \to I^{\vee \vee} \to 
\bigoplus_{i \in J} E_0 \otimes \Hom(E_0,{\Bbb C}_{x_i}) \to 0
\end{equation}
where $F'$ fits in an exact sequence
\begin{equation}
 0 \to F \to F' \to \bigoplus_{i \in K} {\cal E}_{|\{x_i \} \times Y} \to 0.
\end{equation}
Since $\Ext^1(E_0, {\cal E}_{|\{x_i \} \times Y})=0$,
we get that $I^{\vee \vee} \cong F' \oplus 
\bigoplus_{i \in J} E_0 \otimes \Hom(E_0,{\Bbb C}_{x_i})$.
Then we see that $I \cong F' \oplus 
\bigoplus_{i \in J}{\cal E}_{|\{x_i \} \times Y}$, which is a contradiction.
Therefore $I$ is locally free. In the same way as above, we get the relation
\eqref{eq:JHF}. Hence we also get our claim. 
\end{proof}

\begin{lem}\label{lem:bogomolov}
Let $F$ be a $G_2$-twisted semi-stable sheaf with
$v(F)=w:=a w_0+l \varrho_Y-
(d\widehat{H}+\widehat{D}+
(d\widehat{H}+\widehat{D},\widetilde{\xi}_0)\varrho_Y /r_0)$.
Then $\langle w^2 \rangle \geq -2 g^2$, where 
$g:=\gcd(ar_0,lr_0)$.
\end{lem}

\begin{proof}
Let 
\begin{equation}
0 \subset F_1 \subset F_2 \subset \dots \subset F_s=F
\end{equation}
be the Jordan-H\"{o}lder filtration of $F$ with respect to the
$G_2$-twisted stability.
We set $v(F_i/F_{i-1}):=
a_i w_0+l_i \varrho_Y-
(d_i\widehat{H}+\widehat{D}_i+
(d_i\widehat{H}+\widehat{D}_i,\widetilde{\xi}_0)\varrho_Y /r_0)$.
Then we can write $(a_i r_0,l_i r_0)=k_i(ar_0,lr_0)/g$,
where $k_i \in {\Bbb Z}$. Since $\sum_{i=1}^s k_i/g=1$,
we get that $s \leq \sum_{i=1}^s k_i=g$,
which implies that $\langle w^2 \rangle=\sum_{i,j=1}^s \langle
v(F_i/F_{i-1}),v(F_j/F_{j-1}) \rangle \geq \sum_{i,j=1}^s 
(-2)=-2s^2 \geq -2g^2$.
\end{proof}

\begin{prop}\label{prop:IT_0}
We set $v:=lv_0^{\vee}+a \varrho_X+
(dH+D)-(dH+D,\xi_0)\varrho_X /r_0$,
where $D \in \NS(X) \otimes {\Bbb Q} \cap H^{\perp}$.
Let $E$ be a $\mu$-stable sheaf with $v(E)=v$.
If 
\begin{equation}
\frac{d}{r_0 l}>\max \{4lr_0^2+1/(2n r_0 l),
2r_0(\langle v^2 \rangle-(D^2))\},
\end{equation}
then $\IT_0$ holds for $E$ with respect to ${{\cal F}}_{\cal E}$.
\end{prop}

\begin{proof}
Assume that $H^1(X,{\cal E}_{|X \times \{y \}} \otimes E)
=\Ext^1({\cal E}_{|X \times \{y \}} \otimes E,{\cal O}_X)^{\vee} \ne 0$.

(I) We first treat the case where ${\cal E}$ is locally free.
We consider a non-trivial extension
\begin{equation}
 0 \to {\cal E}_{|X \times \{y \}}^{\vee} \to I \to E \to 0.
\end{equation}
Let 
\begin{equation}
0 \subset F_1(I) \subset F_2(I) \subset
\dots \subset F_s(I)=I
\end{equation}
be the Harder-Narasimhan filtration
of $I$ with respect to the $\mu$-semi-stability 
(if $s=1$, then $I$ is $\mu$-semi-stable).
We set $I_i:=F_i(I)/F_{i-1}(I)$ and
$v(I_i):=l_i v_0^{\vee}+a_i \varrho_X-
(d_i{H}+{D}_i+
(d_i{H}+{D}_i,{\xi}_0)\varrho_X /r_0)$.
Then $d_1/l_1 >d_2/l_2> \dots >d_s/l_s$.
In the same way as in the proof of Proposition
\ref{prop:WIT_2}, 
we see that $d_i>0$ and $d/l>d_1/l_1$.
Assume that $s>1$. Then $d_i<d$ for all $i$.
By our assumptions, Lemma \ref{lem:key} implies that $a_i < ad_i/d$.
Then we see that $a=\sum_i a_i<\sum_i ad_i/d=a$,
which is a contradiction.
Thus $s=1$. If $I$ is properly $\mu$-semi-stable, then we can apply 
Lemma \ref{lem:key} again, and we get a contradiction.
If $I$ is $\mu$-stable, then by Lemma \ref{lem:key-rem} (2), 
we get a contradiction.
Therefore we conclude that
$H^1(X,{\cal E}_{|X \times \{y \}} \otimes E)=0$ for all
$y \in Y$.  

(II) We next assume that ${\cal E}$ is not locally free.
We take a locally free resolution 
\begin{equation}
 0 \to V_1 \to V_0 \to {\cal E}_s \to 0
\end{equation}
of ${\cal E}_s$
such that $\Ext^i(V_0,E_0^{\vee})=0$, $i>0$.
Then $\Ext^1(V_1,E_0^{\vee}) \cong \Ext^2({\cal E}_s,E_0^{\vee})
\cong \Hom(E_0^{\vee},{\cal E}_s)^{\vee}=0$.
We may assume that $H^0(X,V_0 \otimes E_0)=0$.
We note that 
\begin{equation}
\begin{split}
\Ext^1({\cal E}_s \otimes E,{\cal O}_X)=&
\Ext^1((V_1 \to V_0) \otimes E,{\cal O}_X)\\
=&\Ext^1(E,V_0^{\vee} \to V_1^{\vee})
\end{split}
\end{equation}
and 
$\Ext^1(E,V_0^{\vee} \to V_1^{\vee})$ parametrizes 
diagrams
\begin{equation}\label{eq:ext-diagram}
\begin{matrix}
& & V_1^{\vee} & &&&&&\\
& & \uparrow & \nwarrow & &&&&\\
0& \to & V_0^{\vee} & \to & I &\to E & \to 0.
\end{matrix}
\end{equation}
Moreover $\phi:I=V_0^{\vee} \oplus E \to V_1^{\vee}$
with $\phi(E)=0$ defines the $0 \in \Ext^1(E,V_0^{\vee} \to V_1^{\vee})$. 
For a diagram \eqref{eq:ext-diagram}, 
we get the following exact and commutative diagram:
\begin{equation}
\begin{CD}
@. 0 @.@.@.\\
@. @AAA @.@.@.@.\\
0 @>>> \im V_0^{\vee} @>>> V_1^{\vee} @>>> {\Bbb C}_s @>>> 0\\
@. @AAA @AAA @AAA @.\\
0 @>>> V_0^{\vee} @>>> I @>>> E @>>> 0\\
@. @AAA @AAA @AAA @.\\
0 @>>> E_0^{\oplus t_0} @>>> I' @>>> E' @>>> 0\\ 
@. @AAA @AAA @AAA @.\\
@.0@.0@.0@.
\end{CD}
\end{equation}
where $I':=\ker(I \to V_1^{\vee})$ and 
$E':=\ker(E \to {\Bbb C}_s)$.
If $E \to {\Bbb C}_s$ is a zero-map, then
$I'$ gives an extension of $E$ by $E_0^{\oplus t_0}$.
By the same argument as case (I), we see that $\Ext^1(E,E_0)=0$.
Hence we get a splitting $E \to I'$, which implies that
\eqref{eq:ext-diagram} is the trivial class. 
If $E \to {\Bbb C}_s$ is non-trivial, then
$I \to V_1^{\vee}$ is surjective.
\begin{claim}\label{claim:split}
$\Hom(I',E_0)=0$.
\end{claim}
Proof of Claim \ref{claim:split}:
We note that
$\Ext^1(V_1^{\vee},E_0)=\Ext^1(E_0^{\vee},V_1) \cong
\Ext^1(V_1,E_0^{\vee})^{\vee}=0$.
Hence $\Hom(I,E_0) \to \Hom(I',E_0)$ is surjective.
On the other hand, by the stability condition on $E$, we see that
$\Hom(I,E_0) \to \Hom(V_0^{\vee},E_0)$ is injective. 
Since $\Hom(V_0^{\vee},E_0)=H^0(X,V_0 \otimes E_0)=0$,
we conclude that $\Hom(I',E_0)=0$.

Since $v(I')=v(E)-\varrho_X+t_0 v(E_0)
=v(E)+v({\cal E}_s)^{\vee}$, applying the same arguments as case (I)
to $I'$,
we see that $\Ext^1(E,V_0^{\vee}\to V_1^{\vee})=0$.
\end{proof}

\begin{cor}\label{cor:IT_0}
Keep notation as above.
Let $E$ be a $G_1$-twisted stable sheaf with $v(E)=v$.
\begin{enumerate}
\item[(1)]
If 
\begin{equation}
\frac{d}{r_0 l}>\max \{4lr_0^2+1,2r_0(\langle v^2 \rangle-(D^2)+(r_0 l)^2/2)\},
\end{equation}
then $\IT_0$ holds for $E$ with respect to ${{\cal F}}_{\cal E}$.
\item[(2)]
If $\langle v^2 \rangle>0$ and $H$ is general with respect to $v$, then 
for a stable sheaf $E$ with 
\begin{equation}
\frac{d}{r_0 l}>\max \{4lr_0^2+1,
2r_0(\langle v^2 \rangle-(D^2))\},
\end{equation}
$\IT_0$ holds with respect to ${{\cal F}}_{\cal E}$.
\end{enumerate}
\end{cor}

\begin{proof}
Let 
\begin{equation}
0 \subset F_1 \subset F_2 \subset \dots \subset F_s=E
\end{equation}
be the Jordan-H\"{o}lder filtration of $E$
with respect to the $\mu$-stability.
We set $E_i:=F_i/F_{i-1}$ and
$v(E_i):=
l_i v_0^{\vee}+a_i \varrho_X+
(d_i H+D_i)-(d_i H+D_i,\xi_0)\varrho_X /r_0$,
where $D_i \in \NS(X) \otimes {\Bbb Q} \cap H^{\perp}$.
We first prove (1).
Since $\langle v(E_j)^2 \rangle \geq -2 \geq -2(\rk E_j)^2$,
by using the equality
$(\langle v^2 \rangle-(D^2))/\rk v=\sum_j
(\langle v(E_j)^2 \rangle-(D_j^2))/\rk E_j$,
we see that $\langle v(E_i)^2 \rangle-(D_i^2)  \leq
\langle v^2 \rangle-(D^2)+2 \rk E_i (\rk E-\rk E_i)$.
Since $\rk E_i(\rk E-\rk E_i) \leq (\rk E)^2/2=(r_0 l)^2/2$, 
we get our claim. 

We next prove (2).
Since $H$ is general with respect to $v$,
$c_1(E)/\rk(E)=c_1(E_i)/\rk E_i$ for all $i$.
We shall show that
$\langle v(E_i)^2 \rangle \leq \langle v^2 \rangle$.
Then our claim follows from Proposition \ref{prop:IT_0}.

If there is not a $\mu$-stable sheaf $G$ such that
$\langle v(G)^2 \rangle=-2$ and $c_1(E)/\rk E=c_1(G)/\rk G$, then 
$\langle v(E_i)^2 \rangle \geq 0$ for all $i$.
Hence $\langle v(E_i)^2 \rangle \leq \langle v^2 \rangle$.
We assume that there is a $\mu$-stable sheaf $G$ such that
$\langle v(G)^2 \rangle=-2$ and $c_1(E)/\rk E=c_1(G)/\rk G$.
It is sufficient to prove the following claim.

\begin{claim}\label{claim:cor:IT_0}
Let $E$ be a $\mu$-semi-stable sheaf such that $c_1(E)/\rk E=c_1(G)/\rk G$ 
and $\Hom(G,E)=0$.
Let 
$0 \subset F_1 \subset F_2 \subset \dots \subset F_s=E$
be the Jordan-H\"{o}lder filtration of $E$
with respect to the $\mu$-stability and set $E_i:=F_i/F_{i-1}$.
Then $\langle v(E_i)^2 \rangle \leq \langle v(E)^2 \rangle$.
\end{claim}
Proof of Claim \ref{claim:cor:IT_0}:
We note that $(\rk G,c_1(G))$ is primitive.
Hence we can set
$v(E)=n v(G)-a \varrho_X$ where $n$ and $a$ are positive integers.
We shall prove our claim by induction on $n$.
Since $E_i$ are $\mu$-stable,
we can write $v(E_i)=n_i v(G)-a_i \varrho_X$,
$a_i \geq 0$. Hence $0 \leq a_i \leq a$ for all $i$.
Then $\langle v(E_i)^2 \rangle \leq 2n_i(a \rk G-n_i)$.
Since $x(a\rk G-x)$ is increasing for $0 \leq x \leq a \rk G/2$,
we see that $\langle v(E_i)^2 \rangle
\leq \langle v(E)^2 \rangle$, if $n_i=1$.
If $\langle v(E),v(G) \rangle=a \rk G-2n \geq 0$,
we also see that $\langle v(E_i)^2 \rangle
\leq \langle v(E)^2 \rangle$ for all $i$.
We assume that $\langle v(E),v(G) \rangle=a \rk G-2n<0$.
Then $k:=\dim \Hom(E,G) \geq -\langle v(E),v(G) \rangle>0$.
We note that
$\phi:E \to G \otimes \Hom(E,G)^{\vee}$ is generically surjective.
Hence we can set that 
$v(\ker \phi):=v(E)-k v(G)+b \varrho_X$, $0 \leq b \leq a$.
Since $(\im \phi)^{\vee \vee}=G^{\oplus k}$ and our claim holds 
for $E_i$ with $n_i=1$, it is sufficient to show our claim
for $E_i$ in $\ker \phi$.
Since
$\langle v(\ker \phi)^2 \rangle=
\langle v(E)^2 \rangle-2k(\langle v(E),v(G) \rangle+k)-b(n-k)\rk G \leq 
\langle v(E)^2 \rangle$ and $\Hom(G,\ker \phi)=0$,
by using the induction hypothesis, we get our claim.
\end{proof}

\begin{rem}\label{rem:IT_0}
If $\NS(X)={\Bbb Z}H$, then the same assertion holds
for an isotropic Mukai vector:
In this case, we may assume that there is a $\mu$-stable vector bundle
$G$ with $\langle v(G)^2 \rangle=-2$ such that
$E=\ker(\Hom(G,{\Bbb C}_x) \otimes G \to {\Bbb C}_x)$.
Then $E$ fits in an exact sequence
\begin{equation}
 0 \to G' \to E \to G^{\oplus (\rk G-1)} \to 0
\end{equation}
where $G'$ is the kernel of a surjective homomorphism
$\psi:G \to {\Bbb C}_x$.
Then we get $\langle v(G')^2 \rangle=2(\rk G-1)$.
Hence $(4 r_0 \rk E+1)-2 r_0 \langle v(G')^2 \rangle=
4r_0(\rk G^2 -\rk G)+4 r_0+1>0$.
Applying Proposition \ref{prop:IT_0} to $G'$,
we see that $\IT_0$ holds for $G'$, and hence for $E$. 
\end{rem}

\begin{prop}\label{prop:IT_0-2}
We set $v:=a \varrho_X+
(dH+D)-(dH+D,\xi_0)\varrho_X /r_0 \in H^{ev}(X,{\Bbb Z})$,
where $a>0$ and $D \in \NS(X) \otimes {\Bbb Q} \cap H^{\perp}$.
Let $E$ be a $G_1$-twisted stable sheaf with $v(E)=v$.
If 
\begin{equation}
a>\max \{2r_0+1,(\langle v^2 \rangle-(D^2))/2+1\},
\end{equation}
then $\IT_0$ holds for $E$ with respect to ${{\cal F}}_{\cal E}$.
\end{prop}

\begin{proof}
Assume that $H^1(X,{\cal E}_{|X \times \{y \}} \otimes E)
=\Ext^1({\cal E}_{|X \times \{y \}} \otimes E,{\cal O}_X)^{\vee} \ne 0$.
We only treat the case where ${\cal E}$ is locally free.
The other case is similar to the proof of Proposition \ref{prop:IT_0}.
We consider a non-trivial extension
\begin{equation}
 0 \to {\cal E}_{|X \times \{y \}}^{\vee} \to I \to E \to 0.
\end{equation}
Assume that $I$ is not $\mu$-semi-stable.
Let 
$I_0$ be the torsion submodule of $I$ and
\begin{equation}
0 \subset F_1(I/I_0) \subset F_2(I/I_0) \subset
\dots \subset F_s(I/I_0)=I/I_0
\end{equation}
the Harder-Narasimhan filtration
of $I/I_0$. 
We set $I_i:=F_i(I/I_0)/F_{i-1}(I/I_0)$ and
$v(I_i):=l_i v_0^{\vee}+a_i \varrho_X+
(d_i{H}+{D}_i-(d_i{H}+{D}_i,{\xi}_0)\varrho_X /r_0)$,
$1 \leq i \leq s$.
Then we see that $d_1/l_1>d_2/l_2> \dots >d_s/l_s>0$.
We also set $v(I_0):=a_0 \varrho_X+(d_0{H}+{D}_0-
(d_0{H}+{D}_0,{\xi}_0)\varrho_X /r_0)$.
We note that the natural homomorphism $I_0 \to I \to E$ is injective.
Since $E$ is $G_1$-twisted stable, $a_0/d_0<a/d$ or $I_0=0$.
If $d_0=d$, then we see that the exact sequence splits.
Hence we get $d_0<d$.
Since $l_i>0$, we get that
$a_i<d_i a/d$ for $i>0$.
Then we see that $a=\sum_{i=0}^s a_i<\sum_{i=0}^s ad_i/d=a$,
which is a contradiction.
Thus $I$ is $\mu$-semi-stable.
If $I$ is properly $\mu$-semi-stable, we also get a contradiction.
Therefore $I$ is $\mu$-stable.
By Lemma \ref{lem:key-rem}, we get a contradiction.
Therefore $H^1(X,{\cal E}_{|X \times \{y \}} \otimes E)=0$ for all
$y \in Y$.  
\end{proof}

\subsection{Asymptotic stability theorem}
\label{subsect:asymptotic}

\begin{prop}\label{prop:rk=0}
Assume that conditions $(\#1,2)$ holds.
Let $E$ be a $G_1$-twisted stable sheaf with
$v(E):=v=a \varrho_X+(dH+D)-(dH+D,\xi_0)\varrho_X /r_0$.
If 
$a > \max \{2r_0+1,(\langle v^2 \rangle-(D^2))/2+1\}$, 
then ${\cal F}_{\cal E}^0(E)$ is $G_2$-twisted stable.
In particular
${\cal F}_{\cal E}$ induces an isomorphism
\begin{equation}
{\cal M}_H^{G_1}(v)^{s} \to 
{\cal M}_{\widehat{H}}^{G_2}({\cal F}_{\cal E}(v))^{s},
\end{equation}
if ${\cal M}_H^{G_1}(v)^{s} \ne \emptyset$.
\end{prop}

\begin{proof}
By Proposition \ref{prop:IT_0-2},
$\IT_0$ holds for $E$.
We assume that
${\cal F}_{\cal E}^0(E)$ is not $G_2$-twisted semi-stable.
Then there is an exact sequence
\begin{equation}
0 \to F_1 \to {\cal F}_{\cal E}^0(E) \to F_2 \to 0
\end{equation}
such that $F_1 (\ne 0)$ is a $G_2$-twisted stable sheaf with
\begin{equation}
v(F_1)=a_1 w_0+l_1 \varrho_Y-
((d_1 \widehat{H}+\widehat{D}_1)+
(d_1\widehat{H}+\widehat{D}_1,\widetilde{\xi}_0)\varrho_Y /r_0),
d_1/a_1 \leq d/a
\end{equation}
 and $F_2 (\ne 0)$ is a torsion free sheaf with $\mu_{\max,G_2}(F_2)<0$.
Applying Lemma \ref{lem:key0},
we see that $l_1 \leq 0$.
Since $F_1$ satisfies $\WIT_2$, we conclude that $l_1=0$.
Since ${\cal F}_{\cal E}^0(E)$ is not $G_2$-twisted semi-stable, 
we may assume that
$d_1/a_1 < d/a$.
Since $\widehat{{\cal F}}_{\cal E}^2(F_1)$ is a torsion sheaf,
$\widehat{{\cal F}}_{\cal E}^1(F_2)$ is also a torsion sheaf.
By Lemma \ref{lem:torsion-free},
$\widehat{{\cal F}}_{\cal E}^1(F_2)=0$.
Then $\widehat{{\cal F}}_{\cal E}^2(F_1)$ is a destabilizing subsheaf of
$E$, which is a contradiction.
Thus ${\cal F}_{\cal E}^0(E)$ is $G_2$-twisted semi-stable.
By the same proof, we also see that
(1) ${\cal F}_{\cal E}^0(E)$ is $G_2$-twisted stable,
provided that $E$ is $G_1$-twisted stable, and (2)
${\cal F}_{\cal E}$ preserves $S$-equivalence classes. 
Hence we have a morphism
$f:\overline{M}_H^{G_1}(v) \to 
\overline{M}_{\widehat{H}}^{G_2}({\cal F}_{\cal E}(v))$.
Let $\overline{M_H^{G_1}(v)}$ (resp. 
$\overline{M_{\widehat{H}}^{G_2}({\cal F}_{\cal E}(v))}$)
be the closure of $M_H^{G_1}(v)$ in $\overline{M}_H^{G_1}(v)$
(resp. $M_{\widehat{H}}^{G_2}({\cal F}_{\cal E}(v))$ in
$\overline{M}_{\widehat{H}}^{G_2}({\cal F}_{\cal E}(v))$).
Then $f$ induces a morphism 
$f':\overline{M_H^{G_1}(v)} \to  
\overline{M_{\widehat{H}}^{G_2}({\cal F}_{\cal E}(v))}$. 
By Corollary \ref{cor:slope-stable},
$M_{\widehat{H}}^{G_2}({\cal F}_{\cal E}(v))$ 
consists of $\mu$-stable sheaves.
Let $H'$ be a general ample divisor on $Y$ such that 
${\Bbb Q}_+ H'$ is very close to ${\Bbb Q}_+ \widehat{H}$.
Then $M_{\widehat{H}}^{G_2}({\cal F}_{\cal E}(v))$
is contained in 
$M_{H'}^{G_2}({\cal F}_{\cal E}(v))=M_{H'}({\cal F}_{\cal E}(v))$.
By the irreducibility of 
$M_{H'}({\cal F}_{\cal E}(v))$
\cite{Y:9},
$M_{\widehat{H}}^{G_2}({\cal F}_{\cal E}(v))$ is also irreducible.
Hence $f'$ is surjective.
Therefore $M_H^{G_1}(v) \to M_{\widehat{H}}^{G_2}({\cal F}_{\cal E}(v))$
is also surjective. Since this morphism is an immersion,
it is an isomorphism.
\end{proof}

\begin{defn}
Let $v$ be a Mukai vector with $\rk v=0$.
A polarization $H$ is general with respect to $v$ and 
$G \in K(X) \otimes {\Bbb Q}$, if
for a $G$-twisted semi-stable sheaf $E$ with $v(E)=v$
and a non-trivial subsheaf $F$ of $E$,
\begin{equation}
\frac{\chi_G(F)}{(c_1(F),H)}=\frac{\chi_G(E)}{(c_1(E),H)}
\text{ if and only if }
v(F) \in {\Bbb Q}v.
\end{equation}
\end{defn}
If $\langle v(G),v \rangle \ne 0$, then there is a general polarization:
For an effective divisor class $\xi \in \NS(X)$, we set
\begin{equation}
D_{\xi}:=\{\xi_1 \in \NS(X)|\text{ $\xi_1$ and $\xi-\xi_1$ are represented as effective divisors and $(\xi_1^2) \geq -2 \epsilon$} \}.
\end{equation}
Then $D_{\xi}$ is a finite set. 
We set $\xi=c_1(v)$.
\begin{itemize}
\item[$(*1)$]
Assume that $(\langle v(G),v \rangle \xi_1-b \xi,H) \ne 0$ 
for all $\xi_1 \in D_{\xi}$
and $b \in {\Bbb Z}$ 
with $0 \leq |b|<| \langle v(G),v \rangle|$ and 
$ \langle v(G),v \rangle \xi_1-b \xi \ne 0$.
\end{itemize}
Then $H$ is a general polarization with respect to $v$ and
$G$. 
   
Assume that $H$ satisfies this condition for $v$ and $G_1$.
Then $H$ also satisfies this condition for $v \exp(mH)$ and $G_1$.
We assume that $a:=-\langle v_0^{\vee},v \rangle/r_0 \gg 
d=\deg_{G_1}(v)/(r_0(H^2))$.
\begin{claim}\label{claim:general}
$\widehat{H}$ is a general polarization
with respect to 
${\cal F}_{\cal E}(v)=a w_0-
(\widehat{\xi}+\frac{(\widehat{\xi},\widetilde{\xi}_0)}{r_0} \rho_Y)$
(cf. Definition \ref{defn:general1}).
\end{claim}
Proof of Claim \ref{claim:general}: Assume that there is a filtration
\begin{equation}
0 \subset F_1 \subset F_2 \subset \dots \subset F_s=F
\end{equation}
 such that
$F_i/F_{i-1}$, $1 \leq i \leq s$ is a $G_2$-twisted stable sheaf with
$v(F_i/F_{i-1})=a_i w_0-
(\widehat{\xi}_i+\frac{(\widehat{\xi}_i,\widetilde{\xi}_0)}{r_0} \rho_Y)$.
Then $\langle v(F_i/F_{i-1})^2 \rangle \geq -2 \epsilon$ for all $i$.
Hence $(\xi_i^2) \geq -2 \epsilon$.
Since $(\xi_i,H)/a_i=(\xi,H)/a >0$, $\xi_i$ is represented by an effective 
divisor.
In particular $\xi_1$ and $\xi-\xi_1$ are effective.
Thus $\xi_1$ belongs to $D_{\xi}$. 
Hence we get our claim.
\qed

Under the assumption $(*1)$ on $H$ and $a$, we see that
${\cal F}_{\cal E}$ induces an isomorphism
\begin{equation}
{\cal M}_H^{G_1}(v)^{ss} \to 
{\cal M}_{\widehat{H}}^{G_2}({\cal F}_{\cal E}(v))^{ss},
\end{equation}
if ${\cal M}_H^{G_1}(v)^{ss} \ne \emptyset$.
For the non-emptyness of ${\cal M}_H^{G_1}(v)^{ss}$,
see Remark \ref{rem:exist}.

The following corollary is a supplement to \cite[Thm. 8.1]{Y:7} and
\cite{Y:9}.
\begin{cor}
Let $X$ be a K3 surface or an abelian surface.
Assume that $\rk v=0$ and $\overline{M}_H(v) \ne \emptyset$. Then
$\overline{M}_H(v)$ is a normal variety, if
$H$ is general with respect to $v$.
Moreover if $X$ is a K3 surface and $v$ is primitive, 
then $\overline{M}_H(v)$ is an irreducible symplectic manifold
which is deformation equivalent
to $\Hilb_X^{\langle v^2 \rangle/2+1}$.
\end{cor}
   
\begin{proof}
If $X$ is an abelian surface, we assume that ${\cal E}$ is the Poincar\'{e}
line bundle on $X \times \widehat{X}$ and
if $X$ is a K3 surface, we assume that ${\cal E}=I_{\Delta}$,
 where $\Delta \subset X \times X$ is the diagonal.
We set $v=\xi+a \varrho_X$.
We assume that $\overline{M}_H(v) \ne \emptyset$.
Since $H$ is general, $\overline{M}_H(v)$ is normal and
$M_H(v)$ is an open dense subscheme of $\overline{M}_H(v)$.
Hence we shall show that $M_H(v)$ is irreducible.
Replacing $v$ by $v \ch(H^{\otimes m})$,
we may assume that $a \gg d=(\xi,H)$.
By Proposition \ref{prop:rk=0},
we have an isomorphism
${M}_H(\xi+a\varrho_X) \to 
{M}_{\widehat{H}}(a-\widehat{\xi})$.
Since ${M}_{\widehat{H}}(a-\widehat{\xi})$ consists of $\mu$-stable
vector bundles, ${M}_{\widehat{H}}(a-\widehat{\xi})$ is contained in
$M_{H'}(a-\widehat{\xi})$, where 
$H'$ is a general ample divisor on $Y$ such that 
${\Bbb Q}_+ H'$ is very close to ${\Bbb Q}_+ \widehat{H}$.
By \cite{Y:9}, it is irreducible.
Hence we get our claim.
\end{proof}

\begin{rem}\label{rem:exist}
We note that a torsion free sheaf on an irreducible
and reduced curve is stable. 
Hence if there is an irreducible
and reduced curve $C$ with $C=c_1(v) \in \NS(X)$,
then $M_H(v)$ is not empty.
We first assume that $X$ is an abelian surface.
Hence if $c_1(v)$ is not primitive, then 
$M_H(v)$ is not empty. If $c_1(v)$ is primitive,
then the non-emptyness comes from \cite{Y:7}.
We next assume that $X$ is a K3 surface.
If $c_1(v)$ is nef, then 
there is an irreducible and reduced curve $C$ with 
$C=c_1(v) \in \NS(X)$, unless
$c_1(v)=\sigma+nf$,
where $\sigma$ is a section of
an elliptic surface $\pi:X \to {\Bbb P}^1$
and $f$ a fiber of $\pi$ (\cite{SD:1}). 
On an elliptic surface $\pi:X \to {\Bbb P}^1$,  
it is easy to construct a stable sheaf on
a curve $C$ with $C=\sigma+nf \in \NS(X)$.
Therefore $M_H(v) \ne \emptyset$, provided that $c_1(v)$ is nef. 
\end{rem}

\begin{thm}\label{thm:asymptotic}
Assume that conditions $(\#1,2)$ holds.
Let $E$ be a $G_1$-twisted semi-stable sheaf
with $v(E)=v:=lv_0^{\vee}+a \varrho_X+
(dH+D)-(dH+D,\xi_0)\varrho_X /r_0$,
where $D \in \NS(X) \otimes {\Bbb Q} \cap H^{\perp}$.
If $l r_0=1,2$ and 
\begin{equation}
\frac{d}{r_0 l}>\max \{4lr_0^2+1,2r_0(\langle v^2 \rangle-(D^2)+(r_0 l)^2/2)\},
\end{equation}
then ${\cal F}_{\cal E}^0(E)$ is $G_2$-twisted semi-stable.
In particular, ${\cal F}_{\cal E}$ induces an isomorphism
\begin{equation}
{\cal M}_H^{G_1}(v)^{s} \to 
{\cal M}_{\widehat{H}}^{G_2}({\cal F}_{\cal E}(v))^{s},
\end{equation}
if ${\cal M}_H^{G_1}(v)^{s} \ne \emptyset$.
\end{thm}

\begin{proof}
By Corollary \ref{cor:IT_0}, $E$ satisfies $\IT_0$.
Assume that there is an exact sequence
\begin{equation}\label{eq:destab}
0 \to F_1 \to {\cal F}_{\cal E}^0(E) \to F_2 \to 0
\end{equation}
such that $F_1 (\ne 0)$ is a $G_2$-twisted stable sheaf with
\begin{equation}
v(F_1)=a_1 w_0+l_1 \varrho_Y-(d_1 \widehat{H}+\widehat{D}_1
+(d_1 \widehat{H}+\widehat{D}_1,\widetilde{\xi}_0)\varrho_Y /r_0),
0<d_1/a_1 \leq d/a
\end{equation}
 and $F_2 (\ne 0)$ is a torsion free sheaf with $\mu_{\max,G_2}(F_2)<0$.
Since $0<a_1<a$, we get $d_1 \leq da_1/a<d$.
Applying Lemma \ref{lem:key} to the sheaf $F_1$,
we get that $l_1 \leq ld_1/d$.
In the exact sequence \eqref{eq:FM1},
Lemma \ref{lem:slope-bound} implies that
$\widehat{{\cal F}}_{\cal E}^0(F_2)=0$.
Hence $\WIT_2$ holds for $F_1$, which implies that $l_1 \geq 0$.
By Lemma \ref{lem:torsion-free}, $\widehat{{\cal F}}_{\cal E}^1(F_2)$
is torsion free. Since $E$ is also torsion free, 
$\widehat{\cal F}_{\cal E}^2(F_1)$ is a torsion free sheaf
of rank $l_1r_0<lr_0 \leq 2$.
If $\rk v=1$, then $\widehat{\cal F}_{\cal E}^2(F_1)=0$,
which is a contradiction.
If $\rk v=2$, then 
$\widehat{\cal F}_{\cal E}^2(F_1)$ is a torsion free sheaf of rank 1.
By the $G_1$-twisted semi-stability of $E$ 
and $(l_1/l)d-d_1 \leq 0$, we see that 
(i) $\widehat{\cal F}_{\cal E}^2(F_1) \to E$ is a $0$-map
or (ii) $d_1/l_1=d/l$, $a_1/l_1 \leq a/l$ and 
$\widehat{\cal F}_{\cal E}^2(F_1) \to E$ is injective.
If the case (i) occurs, then
$\widehat{{\cal F}}_{\cal E}^1(F_2) \cong 
\widehat{{\cal F}}_{\cal E}^2(F_1)$.
Hence $\widehat{{\cal F}}_{\cal E}^1(F_2)$ satisfies $\WIT_0$.
In the same way as in the proof of Lemma \ref{lem:slope-bound},
we get a contradiction.
If the case (ii) occurs, then by the inequality $d_1/a_1 \leq d/a$,
we see that $a_1/l_1 \geq a/l$.
Therefore $a_1/l_1=a/l$, which implies that $d_1/a_1=d/a$.

If ${\cal F}_{\cal E}^0(E)$ is not $G_2$-twisted semi-stable, then
by Lemma \ref{lem:slope-bound}, there is an exact sequence \eqref{eq:destab}
with $0<d_1/a_1<d/a$, which is a contradiction.
Moreover if  ${\cal F}_{\cal E}^0(E)$ is not $G_2$-twisted stable, then
we also see that $E$ is not $G_1$-twisted stable. 
The last claim follows from the same argument as in the proof of
Proposition \ref{prop:rk=0}.
\end{proof}

\begin{rem}
If $H$ satisfies the inequality \eqref{eq:general} and $\langle v^2 \rangle>0$,
then ${\cal F}_{\cal E}$ induces an isomorphism
\begin{equation}
{\cal M}_H^{G_1}(v)^{ss} \to 
{\cal M}_{\widehat{H}}^{G_2}({\cal F}_{\cal E}(v))^{ss},
\end{equation}
if $\rk v=1,2$ and 
\begin{equation}
\frac{d}{r_0 l}>\max \{4lr_0^2+1,2r_0(\langle v^2 \rangle-(D^2))\}.
\end{equation}
\end{rem}

\subsubsection{The case where $\NS(X)={\Bbb Z}$}

In the above theorem, the choice of $d$ depends on $\langle v^2 \rangle$
and $(D^2)$.
Hence if $\NS(X)={\Bbb Z}H$, then the choice depends only on 
$\langle v^2 \rangle$.
Under this assumption, we can show the asymptotic stability 
generally. 

\begin{prop}\label{thm:asymptotic3}
Assume that $\NS(Y)={\Bbb Z}\widehat{H}$.
We set $w:=a w_0+l \varrho_Y-
(d\widehat{H}+
(d\widehat{H},\widetilde{\xi}_0)\varrho_Y /r_0)$.
Let $F$ be a stable sheaf with 
$v(F)=w$.
If 
\begin{equation}\label{eq:cond3}
\frac{d}{r_0 l}>\max \{4lr_0^2+1,2r_0\langle w^2 \rangle\},
\end{equation}
then $\widehat{{\cal F}}_{\cal E}^2(F)$ is stable.
\end{prop}

\begin{proof}
By Proposition \ref{prop:WIT_2}, $\WIT_2$ holds for $F$ and
$\widehat{{\cal F}}_{\cal E}^2(F)$ is torsion free.
Assume that $E:=\widehat{{\cal F}}_{\cal E}^2(F)$ is not semi-stable.
Let 
\begin{equation}
0 \subset F_1(E) \subset F_2(E) \subset \cdots \subset F_s(E)=E
\end{equation}
be the Harder-Narasimhan filtration
of $E$ with respect to semi-stability.
We set $E_i:=F_i(E)/F_{i-1}(E)$ and $v(E_i):=l_i v_0^{\vee}+a_i \varrho_X+
(d_i H-(d_i H,\xi_0)\varrho_X /r_0)$.
By Lemma \ref{lem:slope-bound},
 $d_1/l_1 \geq d_2/l_2 \geq \dots \geq d_s/l_s>0$.
Assume that $d_j/l_j \geq d/l$ for $1 \leq j \leq t$ and
$d_j/l_j<d/l$ for $j>t$.
\begin{claim}\label{claim:t=s}
 $t=s$, that is, $E$ is $\mu$-semi-stable. 
\end{claim}
Proof of Claim \ref{claim:t=s}:
Since $\Hom(E_i,E_j)=0$ for $i<j$, 
\cite[Cor. 2.8]{Mu:4} implies that
\begin{equation}
\sum_{i=1}^s \dim \Ext^1(E_i,E_i) \leq 
\dim \Ext^1(E,E) = \langle w^2 \rangle+2.
\end{equation}
Then $\langle v(E_i)^2 \rangle \leq \dim \Ext^1(E_i,E_i)-2
\leq \langle w^2 \rangle$ for all $i$.
For $E_i$ with $i \leq t$,
we take the Jordan-H\"{o}lder filtration of $E_i$:
\begin{equation}
0 \subset F_1^J(E_i) \subset F_2^J(E_i) \subset \dots \subset
F_{s_i}^J(E_i)=E_i.
\end{equation}
Since $H$ is general with respect to all Mukai vectors, 
$\langle v(F_j^J(E_i)/F_{j-1}^J(E_i))^2 \rangle =(r_j/\rk E_i)^2 
\langle v(E_i)^2 \rangle$, where $r_j:=\rk F_j^J(E_i)/F_{j-1}^J(E_i)$.
Hence $\langle v(F_j^J(E_i)/F_{j-1}^J(E_i))^2 \rangle \leq  
\langle v(E_i)^2 \rangle$, or $E_i=G_i^{\oplus n_i}$, where
$G_i$ is a stable vector bundle with
$\langle v(G_i)^2 \rangle=-2$. 
Applying Corollary \ref{cor:IT_0} (or Remark \ref{rem:IT_0})
to each $F_j^J(E_i)/F_{j-1}^J(E_i)$, 
we see that $\IT_0$ holds for $E_i$, $i \leq t$.
Therefore $F_t(E)$ also satisfies $\IT_0$.
Since $E$ satsifies $\IT_0$, $E/F_t(E)$ also satisfies $\IT_0$ and 
we get an exact sequence
\begin{equation}
 0 \to {\cal F}_{\cal E}^0(F_t(E)) \to F \to {\cal F}_{\cal E}^0(E/F_t(E)) 
\to 0.
\end{equation}
We set $v(E/F_t(E))=l' v_0^{\vee}+a'\varrho_X+
(d' H-(d' H,\xi_0)\varrho_X /r_0)$.
Then $v({\cal F}_{\cal E}^0(E/F_t(E)))=a' w_0+l'\varrho_Y-
(d' \widehat{H}+(d' \widehat{H},\widetilde{\xi_0})\varrho_Y /r_0)$.
Hence $a' \geq 0$ and if $a'=0$, then $d' \leq 0$.
On the other hand, by our assumption
 \eqref{eq:cond3} and Lemma \ref{lem:key} (2), we get that
$a_i<ad_i/d$ for $i>t$, which implies that
$a'<ad'/d$. Then we see that $a'>0$ and ${\cal F}_{\cal E}^0(E/F_t(E))$
gives a destabilizing quotient sheaf of $F$.
Therefore $t=s$. 

Then $d_1/l_1=d/l$ and $a_1/l_1>a/l$.
Since $E_1$ satisfies $\IT_0$, the inequality $-d_1/a_1>-d/a$ implies that 
${\cal F}_{\cal E}^0(E_1)$ is a destabilizing subsheaf of $F$.
Therefore $E$ is semi-stable.
If $E$ is not stable, then
$E$ contains a subsheaf $E_1$ with
\begin{equation}
\text{$v(E_1)=l_1 v_0^{\vee}+a_1 \varrho_X+
(d_1 H-(d_1 H,\xi_0)\varrho_X /r_0)$,
$d_1/l_1=d/l$ and $a_1/l_1 \geq a/l$.}
\end{equation}
Then $\IT_0$ holds for $E_1$, $E/E_1$ and we have an exact sequence
\begin{equation}
 0 \to {\cal F}_{\cal E}^0(E_1) \to F \to {\cal F}_{\cal E}^0(E_2) 
 \to 0.
\end{equation}
Since $d_1/l_1 \leq d/l$, we get a contradiction.
Thus $E$ is stable.
\end{proof}

\begin{thm}\label{thm:asymptotic2}
Assume that $\NS(X)={\Bbb Z}H$. 
We set $v:=lv_0^{\vee}+a \varrho_X+
dH-(dH,\xi_0)\varrho_X /r_0$.
Let $E$ be a stable sheaf with $v(E)=v$.
If 
\begin{equation}
\frac{d}{r_0 l}>\max \{4lr_0^2+1,2r_0\langle v^2 \rangle\},
\end{equation}
then ${{\cal F}}_{\cal E}^0(E)$ is stable.
In particular, ${\cal F}_{\cal E}$ induces an isomorphism
\begin{equation}
{\cal M}_H(v)^{ss} \to 
{\cal M}_{\widehat{H}}({\cal F}_{\cal E}(v))^{ss}.
\end{equation} 

\end{thm}

\begin{proof}
Let $F_1 \subset {\cal F}_{\cal E}^0(E)$ 
be the first filter of the Harder-Narasimhan
filtration of ${\cal F}_{\cal E}^0(E)$. 
We set $v(F_1)=a_1 w_0+l_1 \varrho_{Y}-(d_1 \widehat{H}+
(d_1 \widehat{H},\widetilde{\xi_0})\varrho_Y /r_0)$.
Then $0<d_1/a_1<d/a$, or $d_1/a_1=d/a$ and $l_1/a_1>l/a$.
Since $\Hom(F_1,{\cal F}_{\cal E}^0(E)/F_1)=0$, 
\cite[Cor. 2.8]{Mu:4} implies that
$\dim \Ext^1(F_1,F_1) \leq 
\dim \Ext^1({\cal F}_{\cal E}^0(E),{\cal F}_{\cal E}^0(E))$.
Then $\langle v(F_1)^2 \rangle \leq \dim \Ext^1(F_1,F_1)-2
\leq \langle v^2 \rangle$.
Let $F_1'$ be a stable subsheaf of $F_1$ such that
$v(F_1')=a_1' w_0+l_1' \varrho_{Y}-(d_1' \widehat{H}+
(d_1' \widehat{H},\widetilde{\xi_0})\varrho_Y /r_0)$ with
$d_1'/a_1'=d_1/a_1$ and $l_1'/a_1'=l_1/a_1$.
Then $\langle v(F_1')^2 \rangle=(l_1'/l_1)^2
\langle v(F_1)^2 \rangle$.
Since $d_1/l_1 \geq d/l$, $F_1'$ satisfies the condition \eqref{eq:cond3} in
Proposition \ref{thm:asymptotic3},
and hence $\widehat{\cal F}_{\cal E}^2(F_1')$ is a stable sheaf.
Then by the same argument as in the proof of Theorem \ref{thm:asymptotic},
we see that the claim holds.
\end{proof}

\subsection{A special case} 
Let $(X,H)$ be a polarized abelian surface with $\NS(X)={\Bbb Z}H$.
We set $(H^2)=2n$. 
Let ${\cal P}$ be the Poincar\'{e} line bundle on
$X \times \widehat{X} $.
In this special case, we shall give more precise results.
We first treat positive rank cases.

\subsubsection{Positive rank cases}
\begin{prop}\label{prop:key-special}
For positive integers $r,d,a$,
We set $v:=r+dH+a \varrho_X$ and $k:=\gcd(r,d)>0$.
We take a pair of integers $(r',d')$ such that
$rd'-r'd=-k$ and $0 \leq r'<r$.
If $dn>\max\{\frac{1}{2k}(r'+\frac{(k-1)}{k}r)\langle v^2 \rangle,
\frac{1}{2}\langle v^2 \rangle\}$, then the following assertions hold:
\begin{enumerate}
\item[(1)]
For any stable sheaf $F_1$ with
$v(F_1)=a_1 \pm d_1 \widehat{H}+r_1 \varrho_{\widehat{X}}$,
$0<d_1<d$ and $d_1/a_1 \leq d/a$, we have $r_1 \leq rd_1/d$.
\item[(2)]
For any stable sheaf $E_1$ with $v(E_1)=r_1+d_1 H+a_1 \varrho_X$,
$0<d_1<d$ and $d_1/r_1<d/r$, we have $a_1<ad_1/d$.
\end{enumerate}
\end{prop}

\begin{proof}
We shall prove the claim (1).
We set $s:=\langle v^2 \rangle/2$.
Let $F_1$ be a stable sheaf with
$v(F_1)=a_1-d_1 \widehat{H}+r_1 \varrho_{\widehat{X}}$,
$0<d_1<d$ and $d_1/a_1 \leq d/a$.
If $r_1 \leq 0$, then obviously our claim holds. 
If $r_1>0$, then we see that
\begin{equation}
0 \leq \langle v(F_1)^2 \rangle \leq \frac{2d_1}{rd}(nd(rd_1-r_1d)+r_1 s).
\end{equation}
If $r_1 \geq r$, then we get a contradiction by the inequality
$dn>s$. Assume that $r_1<r$.
If $rd_1-r_1d<0$, then there is a positive integer $m$ such that
$rd_1-r_1d=-km$. Then $r_1-r'm$ is divisible by $r/k$ and 
$r_1-r'm<r$.
Hence we get $r_1-r'm \leq r-r/k$, which implies that
$nd(rd_1-r_1d)+r_1 s =-mknd+r_1 s \leq -mknd+rs-rs/k+r'm s<0$, 
which is a contradiction.
Therefore $rd_1-r_1d \geq 0$.
\end{proof}

Then we get the following.
\begin{thm}\label{thm:asymptotic-special}
${\cal G}_{\cal P}$ induces an isomorphism
\begin{equation}
{\cal M}_H(r+dH+a\varrho_X)^{ss} \to 
{\cal M}_{\widehat{H}}(a+d \widehat{H}+r\varrho_{\widehat{X}})^{ss}
\end{equation}
if $dn > rs$, where $s:=(d^2n-ra)$.
\end{thm}

\begin{proof}
We note that $rs \geq \frac{1}{k}(r'+\frac{(k-1)}{k}r)s$.
Under our conditions, by a modification of the proof of
Proposition \ref{prop:IT_0},
we see that $\WIT_2$ holds with respect to
${\cal G}_{\cal P}$.
Assume that ${\cal G}_{\cal P}^2(E)$ is not semi-stable.
Then we have an exact sequence
\begin{equation}
 0 \to G_1 \to {\cal G}_{\cal P}^2(E) \to G_2 \to 0,
\end{equation}
where $G_1$ is a torsion free sheaf with 
$\mu_{\min,{\cal O}_{\widehat{X}}}(G_1)>0$ and
$G_2$ is a stable sheaf with
$v(G_2)=a_2+d_2 \widehat{H} +r_2 \varrho_{\widehat{X}}$
such that (i) $0<d_2/a_2<d/a$, or (ii) $d_2/a_2=d/a$ and $r_2/a_2<r/a$. 
Then we see that $\WIT_2$ holds for $G_2$ with respect to
$\widehat{\cal G}_{\cal P}$ and
we have an exact sequence
\begin{equation}
0 \to \widehat{\cal G}_{\cal P}^1(G_1) \to
\widehat{\cal G}_{\cal P}^2(G_2) \to E \to 
\widehat{\cal G}_{\cal P}^2(G_1) \to 0.
\end{equation}
Since $d_2 \leq d a_2/a<d$, Proposition \ref{prop:key-special}
implies that $r_2 \leq r d_2/d$.
By the proof of Theorem \ref{thm:asymptotic2},
we see that $\langle v(G_2)^2 \rangle \leq \langle v(E)^2 \rangle$.
Hence $d_2/r_2 \geq d/r >\langle v(E)^2 \rangle \geq \langle v(G_2)^2 \rangle$.
By the same argument as in Proposition \ref{thm:asymptotic3},
we can show that $\widehat{\cal G}_{\cal P}^2(G_2)$ is a stable sheaf.
Then by the proof of Theorem \ref{thm:asymptotic}, we get our theorem. 
\end{proof}

\begin{rem}
Assume that $r \leq 3$.
Under the notation in Proposition \ref{prop:key-special},
if $k=1$ and $dn>r's$, then 
${\cal G}_{\cal P}$ induces an isomorphism
\begin{equation}
{\cal M}_H(r+dH+a\varrho_X)^{ss} \to 
{\cal M}_{\widehat{H}}(a+d \widehat{H}+r\varrho_{\widehat{X}})^{ss}.
\end{equation}
{\it Proof.}
We use the notation in the proof of Theorem \ref{thm:asymptotic-special}.
If $\rk \widehat{{\cal G}}_{\cal P}(G_1)=1$, then 
Lemma \ref{lem:torsion-free} implies that
$\widehat{{\cal G}}_{\cal P}^1(G_1)$ is a line bundle.
Then by the stability of $E$,
we get $\deg(\widehat{{\cal G}}_{\cal P}^1(G_1))>0$,
which implies that 
$\widehat{{\cal G}}_{\cal P}^1(G_1)$ is an ample line bundle.
Hence $\WIT_2$ holds for $\widehat{{\cal G}}_{\cal P}^1(G_1)$
with respect to ${\cal G}_{\cal P}$.
On the other hand, by using the spectral sequence on 
${\cal G}_{\cal P} \circ \widehat{{\cal G}}_{\cal P}(G_1)$, we see that 
${\cal G}_{\cal P}^2(\widehat{{\cal G}}_{\cal P}^1(G_1))=0$.
Therefore $\rk \widehat{{\cal G}}_{\cal P}(G_1) \ne 1$.
Then the proof is similar to that of Theorem \ref{thm:asymptotic}. 

\end{rem}

\begin{rem}
If $r=1$ and $d \geq 2$, then $\IT_0$ holds with respect to
${\cal F}_{\cal P}$ under 
the assumption $2(d-1)n > s$ (cf. \cite[Thm. 1.1]{T:1}).
\end{rem}

If $dn \leq s$, then ${\cal F}_{\cal P}$ does not always 
preserve the stability.

\begin{lem}\label{lem:counter}
Assume that $d=kr+1$ and $dn \leq s \leq (d^2-(d-1)^2/r)n-2r $.
Then there is a $\mu$-stable sheaf 
$E$ with $v(E)=r+dH+\frac{(d^2n-s)}{r}\varrho_X$
such that $E$ satisfies $\IT_0$ with respect to ${\cal F}_{\cal P}$,
but ${\cal F}_{\cal P}^0(E)$ is not $\mu$-semi stable.
\end{lem}

\begin{proof}
We set $v:=r+dH+\frac{(d^2n-s)}{r}\varrho_X$.
We shall find a member 
$E\in {\cal M}_H(v)^{ss}$ such that
${\cal F}_{\cal P}^0(E)$ is not stable. 

\begin{claim}\label{claim:*}
There is a $\mu$-stable sheaf $E$ with $v(E)=v$
such that $H^0(X, E(-kH)) \ne 0$ and $\IT_0$ holds with respect to
${\cal F}_{\cal P}$.
\end{claim}
We first assume this claim and show that
${\cal F}_{\cal P}^0(E)$ is not stable.
We set $F:=\coker({\cal O}_X \to E(-kH))$.
Then we have an exact sequence
\begin{equation}
0 \to {\cal F}_{\cal P}^0({\cal O}_X(kH)) \to
 {\cal F}_{\cal P}^0(E) \to  
{\cal F}_{\cal P}^0(F(kH)) \to 0.
\end{equation}
Since $v({\cal F}_{\cal P}^0({\cal O}_X(kH))) =
nk^2-k\widehat{H}+\varrho_X$, we get that
\begin{equation}
\begin{split}
\frac{\deg({\cal F}_{\cal P}^0({\cal O}_X(kH)))}
{\rk({\cal F}_{\cal P}^0({\cal O}_X(kH)))} 
-\frac{\deg({\cal F}_{\cal P}^0(E(kH)))}
{\rk({\cal F}_{\cal P}^0(E(kH)))}&=
\frac{-k(H^2)}{k^2n}-\frac{-rd(H^2)}{d^2n-s}\\
&=\frac{2(s-dn)}{k(d^2n-s)} \geq 0.
\end{split}
\end{equation}
Thus ${\cal F}_{\cal P}^0(E)$ is not stable.
Therefore we get our lemma. 

Proof of Claim \ref{claim:*}:
We note that $s \geq n$.
Let $F$ be a stable vector bundle with
$v(F)=(r-1)+{H}-\{(s-n)/r\}\varrho_{X}$.
Then $\Ext^1(F \otimes {\cal P}_{|X \times \{y \}},{\cal O}_X)
=H^1(X,F \otimes {\cal P}_{|X \times \{y \}})^{\vee} \ne 0$
for some $y \in \widehat{X}$.
Let $E$ be a sheaf on $X$ such that $E(-kH)$ is defined as
a non-trivial extension
\begin{equation}
0 \to {\cal O}_X \to E(-kH) \to F \otimes {\cal P}_{|X \times \{y \}} \to 0.
\end{equation}
Then $E$ is $\mu$-stable (see \cite[Lem. 2.1]{Y:5}).
Moreover, since $\chi(F(kH))=(d^2n-s)/r-nk^2=
((d^2-(d-1)^2/r)n-s)/r \geq 2$, Theorem \ref{thm:birat1}
in section \ref{sect:birat} implies that $\IT_0$ holds 
for a general $F$ with respect to ${\cal F}_{\cal P}$.
\end{proof}

\begin{rem}
As we shall see in section \ref{sect:birat},
${\cal F}_{\cal P}$ preserves the stability condition 
for a general $\mu$-stable sheaf.
\end{rem}

%Combining all together, we get the following.
%\begin{cor}
%${\cal F}_{\cal P}$ induces an isomorphism
%\begin{equation}
%{\cal M}_H(1+dH+a\varrho_X)^{ss} \to 
%{\cal M}_{\widehat{H}}(a-d \widehat{H}+\varrho_{\widehat{X}})^{ss}
%\end{equation}
%if and only if $dn > d^2n-a$.
%\end{cor}

\subsubsection{Rank 0 case}

We next treat the rank 0 case.
We start with the following lemma whose proof is similar.
\begin{lem}
We set $v:=d H+a \varrho_X$. If $a>d(d-1)n$, then
\begin{enumerate}
\item[(1)]
for any stable sheaf $F_1$ with
$v(F_1)=a_1 \pm d_1 \widehat{H}+r_1 \varrho_{\widehat{X}}$,
$0<d_1<d$ and $d_1/a_1 \leq d/a$, we have $r_1 \leq 0$, and
\item[(2)]
for any stable sheaf $E_1$ with $v(E_1)=r_1+d_1 H+a_1 \varrho_X$,
$0<d_1<d$ and $r_1>0$, we have $a_1<ad_1/d$.
\end{enumerate}
\end{lem}

\begin{prop}\label{prop:asymptotic4}
${\cal G}_{\cal P}$ induces an isomorphism
\begin{equation}
{\cal M}_H(dH+a\varrho_X)^{ss} \to 
{\cal M}_{\widehat{H}}(a+d \widehat{H})^{ss},
\end{equation}
if $a>d(d-1)n$. Moreover ${\cal F}_{\cal P}$
induces an isomorphism ${\cal M}_H(dH+a\varrho_X)^{ss} \to 
{\cal M}_{\widehat{H}}(a-d \widehat{H})^{ss}$, if $a>d^2n$.
\end{prop}

\begin{proof}
We shall only prove the first claim.
For $E \in {\cal M}_H(dH+a\varrho_X)^{ss}$, we see that $E$ satsifies 
$\WIT_2$ with respect to ${\cal G}_{\cal P}$.
we assume that
${\cal G}_{\cal P}^2(E)$ is not semi-stable.
Then there is an exact sequence
\begin{equation}
0 \to G_1 \to {\cal G}_{\cal P}^2(E) \to G_2 \to 0
\end{equation}
such that 
$G_1$ is a torsion free sheaf with $\mu_{\min,{\cal O}_{\widehat{X}}}(G_1)>0$
and 
$G_2$ is a stable sheaf with
$v(G_2)=a_2+d_2 \widehat{H}+r_2\varrho_{\widehat{X}}$,
(i) $0<d_2/a_2 < d/a$, or (ii) $d_2/a_2=d/a$ and $r_2>0$.
We note that 
\begin{equation}\label{eq:a_1}
a_2 \geq d_2 a/d>d_2(d-1)n.
\end{equation}
If $r_2>0$, then we see that
$0 \leq d_2^2(H^2)-2r_2 a_2 \leq 2nd_2^2-2a_2
<  2nd_2(d_2-d+1) \leq 0$.
Hence we get $r_2 \leq 0$.
Thus the case (ii) does not occur.
Since $\widehat{{\cal G}}_{\cal P}^1(G_1)$ is locally free
(cf. Lemma \ref{lem:torsion-free}) and 
$\widehat{{\cal G}}_{\cal P}^1(G_1)$ 
is a subsheaf of $\widehat{{\cal G}}_{\cal P}^2(G_2)$,
we get $\widehat{{\cal G}}_{\cal P}^1(G_1)=0$. Since $a_2/d_2>a/d$, 
$\widehat{{\cal G}}_{\cal P}^2(G_2)$ is a destabilizing subsheaf of $E$.
Therefore ${\cal G}_{\cal P}^2(E)$ is semi-stable.
\end{proof}

\section{Birational maps}\label{sect:birat}

Let $(X,H)$ be a polarized abelian surface with $\NS(X)={\Bbb Z}H$
again.

\begin{prop}\label{prop:wit-birat}
We set $v:=r+dH+a \varrho_X$, $r,d>0$.
If $\langle v^2 \rangle<2r$, then 
$\WIT_2$ holds for all $\mu$-semi-stable sheaf $E$
with $v(E)=v$.
\end{prop}

\begin{proof}
We shall prove our claim by induction on
$\langle v^2 \rangle$.
Obviously our claim holds for semi-homogenous sheaf.
Let $E$ be a $\mu$-semi-stable sheaf with $v(E)=v$.
Assume that $E$ is $S$-equivalent to $\bigoplus_{i=1}^s E_i$,
where $E_i$, $1 \leq i \leq s$ are $\mu$-stable sheaves.
Then 
\begin{equation}
\sum_i \frac{\langle v(E_i)^2 \rangle}{\rk E_i}=
\frac{\langle v^2 \rangle}{r}<2.
\end{equation}
Since $\langle v(E_j)^2 \rangle \geq 0$ for all $j$,
we get $\frac{\langle v(E_i)^2 \rangle}{\rk E_i}\leq 
\frac{\langle v^2 \rangle}{r}<2$.
Therefore we shall prove our claim for
$\mu$-stable sheaves.

If $a=0$, then $2nd^2<2r$.
Hence the claim follows from Proposition \ref{prop:asymptotic4}.
We assume that $a>0$.
Assume that $\Ext^1(E,{\cal P}_{|X \times \{y\}}) \ne 0$, $y \in \widehat{X}$.
We take a non-trivial extension
\begin{equation}
0 \to {\cal P}_{|X \times \{y\}} \to G \to E \to 0.
\end{equation}
Assume that $G$ is not $\mu$-semi-stable.
Let 
\begin{equation}
0 \subset F_1 \subset F_2 \subset \dots \subset F_s=G
\end{equation}
 be the Harder-Narasimhan filtration of $G$
with respect to the $\mu$-semi-stability.
We set $v(F_i/F_{i-1}):=r_i+d_i H+a_i \varrho_X$.
Then $0<d_s/r_s<\dots<d_2/r_2<d_1/r_1<d/r$ and $r_i \leq r$.
We see that
\begin{equation}
 \begin{split}
  \sum_i \frac{\langle v(F_i/F_{i-1})^2 \rangle}{r_i}&<
  \sum_i 2(n\frac{d}{r}d_i-a_i)\\
  &=2n\frac{d^2}{r}-2a=\frac{\langle v^2 \rangle}{r}<2.
 \end{split}
\end{equation}
Since $\langle v(F_j/F_{j-1})^2 \rangle \geq 0$ for all $j$,
we get $\frac{\langle v(F_i/F_{i-1})^2 \rangle}{r_i}
<\frac{\langle v^2 \rangle}{r}<2$.
Since $r_i \leq r$, we get 
$\langle v(F_i/F_{i-1})^2 \rangle <\langle v^2 \rangle$.
By induction hypothesis, our claim holds for $F_i/F_{i-1}$.
Hence $G$ satisfies $\WIT_2$ with respect to ${\cal G}_{\cal P}$.
Since ${\cal P}_{|X \times \{y\}}$ also satisfies $\WIT_2$ 
with respect to ${\cal G}_{\cal P}$, $E$ satisfies $\WIT_2$ 
with respect to ${\cal G}_{\cal P}$.

Assume that $G$ is $\mu$-semi-stable.
Since $\langle v(G)^2 \rangle
=\langle v^2 \rangle-2a$, by induction hypothesis, our claim holds for $G$.
Therefore $E$ satisfies $\WIT_2$ 
with respect to ${\cal G}_{\cal P}$.
\end{proof}

\begin{lem}
We set $v=r+dH+a \rho_X$, $r,d>0$.
\begin{enumerate}
\item[(1)]
If $a>0$, then there is a stable sheaf $E$ with $v(E)=v$ such that
$\Ext^1(E,{\cal O}_X)=H^1(X,E)^{\vee}=0$.
In particular, $\WIT_2$ holds for $E$ with respect to ${\cal G}_{\cal P}$.
\item[(2)]
If $a \leq 0$, then there is a stable sheaf $E$ with $v(E)=v$ such that
$H^0(X,E)=0$.
In particular, $\WIT_1$ holds for $E$ with respect to ${\cal F}_{\cal P}$.
\end{enumerate}
\end{lem}

\begin{proof}
We take an integer $b$ such that
$0 \leq \langle (r+dH+(a+b) \varrho_X)^2 \rangle=\langle v^2 \rangle-2rb<2r$.
We note that $b \geq 0$.
Let $F$ be a stable sheaf with $v(F)=r+dH+(a+b)\varrho_X$ such that
$H^1(X,F)=0$.
We consider a surjective homomorphism
$\phi:F \to \bigoplus_{i=1}^b {\Bbb C}_{x_i}$, where $x_1,x_2,\dots,x_b \in X$.
If we choose a sufficiently general $\phi$,
then 
\begin{equation}
\begin{cases}
H^1(X,\ker \phi)=0, \text{ if $a>0$},\\
H^0(X,\ker \phi)=0, \text{ if $a \leq 0$}.
\end{cases}
\end{equation}
Since $\ker \phi$ is $\mu$-semi-stable,
by the dimension counting in \cite[sect. 2]{Y:8}, we see that 
a $\mu$-semi-stable sheaf deforms to a stable sheaf.
Hence we get our claim.
\end{proof}

By \cite[Cor. 4.15]{Y:7}, we get the following theorem
which was conjectured in \cite[Conj. 4.16]{Y:7}.
\begin{thm}\label{thm:birat1}
Assume that $r,d>0$. 
\begin{enumerate}
\item[(1)]
If $a>0$, then
${\cal G}_{\cal P}$ induces a birational map
\begin{equation}
\overline{M}_H(r+dH+a \varrho_X) \cdots \to 
\overline{M}_{\widehat{H}}(a+d \widehat{H}+r \varrho_{\widehat{X}}).
\end{equation}
\item[(2)]
If $a \leq 0$, then
${\cal F}_{\cal P}$ induces a birational map
\begin{equation}
\overline{M}_H(r+dH+a \varrho_X) \cdots \to 
\overline{M}_{\widehat{H}}(-a+d \widehat{H}-r \varrho_{\widehat{X}}).
\end{equation}
\end{enumerate}
\end{thm}

\begin{defn}
For a divisor $D$ on $X$,
we define $T_D:{\bf D}(X) \to {\bf D}(X)$ by
sending $F \in {\bf D}(X)$ to $F \otimes {\cal O}_X(D) \in {\bf D}(X)$.
\end{defn}

\begin{lem}\label{lem:birat3}
If $r>0$, then $\overline{M}_H(r+dH+a \varrho_X)$ is birationally
equivalent to $\overline{M}_H(r-dH+a \varrho_X)$.
\end{lem}

\begin{proof}
If $r|d$, then $T_{-2dH/r}$ induces an isomorphism
$\overline{M}_H(r+dH+a \varrho_X) \to \overline{M}_H(r-dH+a \varrho_X)$.
If $r \nmid d$, then there is a $\mu$-stable vector bundle $E$ with
$v(E)=r+dH+a \varrho_X$ (cf. \cite[sect. 2]{Y:8}). 
Since $E^{\vee}$ is also $\mu$-stable,
we get a desired birational map. 
\end{proof}

The following was proved in \cite{Y:7}. 
\begin{prop}\cite[Thm. 9.4]{Y:7}\label{prop:birat2}
If $r, b>0$, then ${\cal G}_{\cal P}$ induces an isomorphism
$\overline{M}_H(r-b \varrho_X) \to 
\overline{M}_{\widehat{H}}(b-r \varrho_{\widehat{X}})$.
\end{prop} 

\begin{defn}
\begin{enumerate}
\item[(1)]
A Mukai vector $v:=r+dH+a \varrho_X$ is positive,
if (i) $r>0$, or (ii) $r=0$ and $d>0$, or
(iii) $r=d=0$ and $a>0$.
We denote a positive $v$ by $v>0$.
\item[(2)]
For a Mukai vector $v$ with $-v>0$,
we set $\overline{M}_H(v):=\overline{M}_H(-v)$.
\end{enumerate}
\end{defn}
By using Theorem \ref{thm:birat1}, Lemma \ref{lem:birat3}
and Proposition \ref{prop:birat2}, we get the following theorem. 
\begin{thm}
For a Mukai vector $v$,
$\overline{M}_H(v)$ is birationally equivalent to
$\overline{M}_{\widehat{H}}({\cal F}_{\cal P}(v))$.
\end{thm}

Assume that $(X,H)$ is a principally polarized abelian surface,
i.e, $n=1$.
We identify $\widehat{X}$ with $X$ by the canonical morphism
$\phi_H:X \to \widehat{X}$. 
Mukai \cite[Thm. 3.13]{Mu:2} showed that $SL(2,{\Bbb Z})$ acts on
${\bf D}(X)$ up to shift such that the correspondence is given by 
\begin{equation}
\begin{matrix}
\begin{pmatrix}
0 & 1\\
-1 & 0
\end{pmatrix} 
\mapsto {\cal F}_{\cal P},  
&
\begin{pmatrix}
1 & 1\\
0 & 1
\end{pmatrix} 
\mapsto T_H.
\end{matrix}
\end{equation} 
Hence we get an $SL(2,{\Bbb Z})$ action on
$H^*(X,{\Bbb Z})_{alg}:={\Bbb Z} \oplus \NS(X) \oplus {\Bbb Z}\varrho_X$.

\begin{cor}
Let $(X,H)$ be a principally polarized abelian surface with
$\NS(X)={\Bbb Z}$.
Then for $g \in SL(2,{\Bbb Z})$ and $v \in H^*(X,{\Bbb Z})_{alg}$,
$\overline{M}_H(v)$ is birationally equivalent to 
$\overline{M}_H(g(v))$.
\end{cor}


\begin{thebibliography}{[Ma-Yk]}
\bibitem[BBH]{BBH:1}
Bartocci, C., Bruzzo, U., Hern\'{a}ndez Ruip\'{e}rez, D.,
{\it A Fourier-Mukai transform for stable bundles on $K3$ surfaces,}
J. Reine Angew. Math. {\bf 486} (1997), 1--16
\bibitem[Br]{Br:2}
Bridgeland, T.,
{\it Equivalences of triangulated categories and Fourier-Mukai
transforms,} 
Bull. London Math. Soc. {\bf 31} (1999), 25--34,
math.AG/9809114
\bibitem[Mu1]{Mu:2}
Mukai, S.,
{\it Duality between $D(X)$ and $D(\hat{X})$ with its application
to Picard sheaves,}
Nagoya Math. J., {\bf 81} (1981), 153--175
\bibitem[Mu2]{Mu:5}
Mukai, S.,
{\it Fourier functor and its application to the moduli of bundles
on an Abelian variety,}
Adv. Studies in Pure Math. {\bf 10} 
(1987), 515--550
\bibitem[Mu3]{Mu:4}
Mukai, S.,
{\it On the moduli space of bundles on K3 surfaces I,}
Vector bundles on Algebraic Varieties, Oxford, 1987, 341--413 
\bibitem[SD]{SD:1}
Saint-Donat, B.,
{\it Projective models of $K3$ surfaces,}
Amer. J. Math. 96 (1974), 602--639
\bibitem[T]{T:1}
Terakawa, H.,
{\it The $k$-very ampleness and $k$-spannedness on polarized abelian
surfaces,}
Math. Nachr. 195 (1998), 237--250

\bibitem[V]{V:1}
Verbitsky, M.,
{\it Projective bundles over hyperkaehler manifolds 
and stability of Fourier-Mukai transform,}
math.AG/0107196 v3

\bibitem[Y1]{Y:5}
Yoshioka, K.,
{\it Some examples of Mukai's reflections on K3 surfaces,}
J. reine angew. Math. {\bf 515} (1999), 97--123   
\bibitem[Y2]{Y:8}
Yoshioka, K.,
{\it Irreducibility of moduli spaces of vector bundles on K3 surfaces,}
math.AG/9907001

\bibitem[Y3]{Y:7}
Yoshioka, K.,
{\it Moduli spaces of stable sheaves on abelian surfaces,}
Math. Ann. {\bf 321} (2001), 817--884, math.AG/0009001


\bibitem[Y4]{Y:9}
Yoshioka, K.,
{\it Twisted stability and Fourier-Mukai transform,}
math.AG/0106118
\end{thebibliography}
\end{document}